\documentclass{article}
\usepackage{amssymb,amsmath,amsthm,tikz,multirow,nccrules,geometry}
\usetikzlibrary{calc, arrows, arrows.meta}
\title{Tiling of Hyperbolic Surface by a Single Tile}
\author{Chunlin Li, Erxiao Wang\thanks{Corresponding author (wang.eric@zjnu.edu.cn).  Research was supported by National Natural Science Foundation of China NSFC-RGC 12361161603 and Key Projects of Zhejiang Natural Science Foundation LZ22A010003.}, Wu Jie,
Zhejiang Normal University \\
Min Yan\thanks{Research was supported by NSFC-RGC Joint Research Scheme N-HKUST607/23 and Hong Kong RGC General Research Fund 16305920.}, 
Hong Kong University of Science and Technology}

\usepackage[hidelinks, hyperindex]{hyperref}

\hyperref[sec:function]{}

\newcommand{\mc}{\mathcal}
\newcommand{\bb}{\mathbb}

\makeatletter
\newsavebox\myboxA
\newsavebox\myboxB
\newlength\mylenA

\newcommand*\xar[2][0.7]{%
    \sbox{\myboxA}{$\m@th#2$}%
    \setbox\myboxB\null
    \ht\myboxB=\ht\myboxA%
    \dp\myboxB=\dp\myboxA%
    \wd\myboxB=#1\wd\myboxA
    \sbox\myboxB{$\m@th\overline{\copy\myboxB}$}
    \setlength\mylenA{\the\wd\myboxA}
    \addtolength\mylenA{-\the\wd\myboxB}%
    \ifdim\wd\myboxB<\wd\myboxA%
       \rlap{\hskip 0.5\mylenA\usebox\myboxB}{\usebox\myboxA}%
    \else
        \hskip -0.5\mylenA\rlap{\usebox\myboxA}{\hskip 0.5\mylenA\usebox\myboxB}%
    \fi}
\makeatother

\makeatletter
\newsavebox\myboxC
\newsavebox\myboxD
\newlength\mylenC

\newcommand*\yar[2][0.7]{%
    \sbox{\myboxC}{$\m@th#2$}%
    \setbox\myboxD\null
    \ht\myboxD=\ht\myboxC%
    \dp\myboxD=\dp\myboxC%
    \wd\myboxD=#1\wd\myboxC
    \sbox\myboxD{$\m@th\underline{\copy\myboxD}$}
    \setlength\mylenC{\the\wd\myboxC}
    \addtolength\mylenC{-\the\wd\myboxD}%
    \ifdim\wd\myboxD<\wd\myboxC%
       \rlap{\hskip 0.5\mylenC\usebox\myboxD}{\usebox\myboxC}%
    \else
        \hskip -0.5\mylenC\rlap{\usebox\myboxC}{\hskip 0.5\mylenC\usebox\myboxD}%
    \fi}
\makeatother

\makeatletter
\newcommand{\drawhypgeodesic}[5][]{
  
  \pgfmathsetmacro{\xone}{#2}
  \pgfmathsetmacro{\yone}{#3}
  \pgfmathsetmacro{\xtwo}{#4}
  \pgfmathsetmacro{\ytwo}{#5}
  \pgfmathsetmacro{\denom}{\xone*\ytwo - \xtwo*\yone}

  \ifdim\denom pt < 0.00001pt
    \ifdim\denom pt > -0.00001pt
      \draw[#1] (#2,#3) -- (#4,#5);
      \def\isLine{1}
    \else
      \def\isLine{0}
    \fi
  \else
    \def\isLine{0}
  \fi
  
  \ifnum\isLine=0

    \pgfmathsetmacro{\mone}{1 + \xone*\xone + \yone*\yone}
    \pgfmathsetmacro{\mtwo}{1 + \xtwo*\xtwo + \ytwo*\ytwo}
    \pgfmathsetmacro{\dd}{(\yone*\mtwo - \ytwo*\mone)/\denom}
    \pgfmathsetmacro{\ee}{(\xtwo*\mone - \xone*\mtwo)/\denom}
    \pgfmathsetmacro{\cx}{-\dd/2}
    \pgfmathsetmacro{\cy}{-\ee/2}
    \pgfmathsetmacro{\rr}{sqrt(\cx*\cx + \cy*\cy - 1)}
    \pgfmathsetmacro{\thetaone}{atan2(\yone-\cy, \xone-\cx)}
    \pgfmathsetmacro{\thetatwo}{atan2(\ytwo-\cy, \xtwo-\cx)}   
    \pgfmathsetmacro{\dtheta}{\thetatwo - \thetaone}
    \pgfmathsetmacro{\dtheta}{\dtheta > 180 ? \dtheta - 360 : \dtheta}
    \pgfmathsetmacro{\dtheta}{\dtheta <= -180 ? \dtheta + 360 : \dtheta}
    
    \draw[#1] (#2, #3) arc (\thetaone : \thetaone + \dtheta : \rr);
  \fi
}
\makeatother

\newtheorem{theorem}{Theorem}
\newtheorem{lemma}[theorem]{Lemma}

\newtheorem{proposition}[theorem]{Proposition}
\newtheorem*{theorem*}{Theorem}

\theoremstyle{definition}
\newtheorem*{definition*}{Definition}
\newtheorem*{case*}{Case}
\newtheorem*{subcase*}{Subcase}

\theoremstyle{remark}

\numberwithin{equation}{section}

\begin{document}

\maketitle

\begin{abstract}
Tilings of a surface of negative Euler characteristic by $n$-gons with $n\ge 7$ is a finite problem. One extreme of the finite problem is single tile tilings. We develop the algorithm for finding all the single tile tilings and present the results for surfaces of small genus.
\end{abstract}

\section{Introduction}
\label{intro}

Consider an edge-to-edge tiling of a closed and connected surface of Euler number $\chi$ by congruent $n$-gons, such that all vertices have degree $\ge 3$. We find that, if $\chi<0$ and $n\ge 7$, then the number of tilings is finite for each fixed surface. Therefore it is possible to develop computer algorithms to find all such tilings. 

The problem has two extreme cases. The first is few number of tiles. The second is maximal number of tiles, which we will see is equivalent to that all vertices have degree 3. In this paper, we develop the algorithm optimised to single tile tilings. Then we apply the algorithm to get all the single tile tilings of the orientable surfaces $2{\bb T}^2$, $3{\bb T}^2$, $4{\bb T}^2$, and the non-orientable surfaces $3{\bb P}^2$, $4{\bb P}^2$, $5{\bb P}^2$. We find the number of tilings quickly becomes very large. For example, the number of tilings of $4{\bb T}^2$ by a single $22$-gon is already 10404513. We can provide the complete list in this paper only for some limited cases. 

A single tile tiling is combinatorially the same as a planar diagram (satisfying degree $\ge 3$ assumption, in this paper). In \cite{zamo1} and \cite{zamo2}, Zamorzaeva-Orleanschi attempted to find all the planar diagrams for $2{\bb T}^2$. Unfortunately, her list is not quite complete, and also includes duplicates. We will compare our findings with hers.

After the combinatorial classification, we should be concerned whether the tilings can be geometrically realised by polygons with straight edges. We believe this is always true for single tile tilings, but cannot come up with a general argument. We explain a way to find examples of such polygons, and explicitly show the geometrical realisation for some families of tilings. 

In subsequent papers \cite{lwwy2} and \cite{lwwy3}, we will study tilings by two tiles, and tilings with maximal number of tiles. Incidentally, not all combinatorial tilings by two tiles can be geometrically realised.

In the rest of the introduction, we develop some basic combinatorial facts. 

Let $v$, $e$, $f$ be the numbers of vertices, edges and tiles in the tiling. Let $v_k$ be the number of vertices of degree $k$. Then we have
\begin{align}
v-e+f &=\chi, \label{eq1} \\
v &= v_3+v_4+v_5+\cdots, \label{eq2} \\
nf=2e &= 3v_3+4v_4+5v_5+\cdots. \label{eq3}
\end{align}
The following is an easy consequence of \eqref{eq3}.

\begin{lemma}\label{parity}
If $n$ is odd, then the number of tiles in a tiling of surface by $n$-gons is even.
\end{lemma}

Canceling $v,e,f$ in \eqref{eq1}, \eqref{eq2}, \eqref{eq3}, we get
\begin{equation}\label{eq4}
-2n\chi
=(n-6)v_3
+(2n-8)v_4
+(3n-10)v_5
+\cdots 
=\sum_{k\ge 3}((k-2)n-2k)v_k. 
\end{equation}
For $k\ge 3$, we have
\[
\frac{n-6}{3}k
=\left(\frac{1}{3}n-2\right)k\le (k-2)n-2k<(n-2)k,
\]
and the inequality $\le$ becomes an equality only for $k=3$. Then by comparing \eqref{eq3} and \eqref{eq4}, we get
\[
\frac{n-6}{3}nf \le -2n\chi< (n-2)nf.
\]
Moreover, the inequality $\le$ becomes an equality only if $v_k=0$ for all $k>3$.

\begin{lemma}\label{bound}
In an edge-to-edge tilings of a surface of Euler number $\chi<0$ by $n$-gons, if $n\ge 7$, then the number $f$ of tiles is bounded
\[
\frac{-2\chi}{n-2}<f\le \frac{-6\chi}{n-6}.
\]
Moreover, $f$ equals the maximal number $\frac{-6\chi}{n-6}$ if and only if all vertices have degree $3$.
\end{lemma}

By Lemma \ref{parity}, we have $f\ge 2$ for odd $n$. Of course we have $f\ge 1$ for even $n$. Combined with $f\le \frac{-6\chi}{n-6}$ in Lemma \ref{bound}, we get
\[
n\le\begin{cases}
3(2-\chi), &\text{odd }n \\
6(1-\chi), &\text{even }n 
\end{cases}.
\]
Combined with the bound for $f$, we know that, for a fixed hyperbolic surface, the study of tilings of the surface by $n$-gons with $n\ge 7$ is a finite problem.

For example, the surface with $\chi=-1$ is $3{\bb P}^2$. We have $n\le 12$. Then for each $7\le n\le 12$, by Lemmas \ref{parity} and \ref{bound}, we get the following complete list of combinations of $f$ and $n$
\begin{align*}
n=10, 12 &\colon f=1; \\
n=9 &\colon f=2; \\
n=8 &\colon f=1,2,3; \\
n=7 &\colon f=2,4,6.
\end{align*}
Moreover, $f$ reaches the respective maximal numbers $6,3,2,1$ for $n=7,8,9,12$. These are exactly the cases when all vertices have degree $3$.

For another example, the surface with $\chi=-2$ is $2{\bb T}^2$ and $4{\bb P}^2$. We similarly get $7\le n\le 18$ and the following complete list
\begin{align*}
n=14,16,18 &\colon f=1; \\
n=12 &\colon f=1,2; \\
n=11 &\colon f=2; \\
n=10 &\colon f=1,2,3; \\
n=9 &\colon f=2,4; \\
n=8 &\colon f=1,2,3,4,5,6; \\
n=7 &\colon f=2,4,6,8,10,12.
\end{align*}
Moreover, $f$ reaches the respective maximal numbers $12,6,4,3,2,1$ for $n=7,8,9,10,12,18$. These are exactly the cases when all vertices have degree $3$. 

The extremes of the problem are the case of minimal $f=1$ or $2$, and the maximal $f=\frac{-6\chi}{n-6}$. In this paper, we describe the computer algorithm for calculating all tilings for $f=1$, i.e., single tile tilings. This implies $v=e-1+\chi\ge 1$, and $n=2e\ge 2(2-\chi)$. Since $n$ is even, we conclude the range for single tile tilings
\[
2(2-\chi)\le n\le 6(1-\chi).
\]
We note the two extreme cases of $n$:
\begin{itemize}
\item $n=2(2-\chi)$ if and only if $v=1$. In other words, the tiling has a single vertex.
\item $n=6(1-\chi)$ if and only if $f=1=\frac{-6\chi}{n-6}$. In other words, all the vertices in the tiling have degree $3$.
\end{itemize}

\section{Single Tile Tiling of Orientable Surface}

A single tile tiling has $f=1$. By \eqref{eq3}, we know $n=2e$ is even. The surface is obtained by glueing the edges of the $n$-gon in $e$ pairs. Therefore a single tile tiling is exactly a {\em planar diagram} $D$ of the $n$-gon. The glueing of the edge pairs gives a single tile tiling ${\mc T}_D$ of a closed surface $S_D$.

We label the {\em corners} of the $n$-gon by $i\in {\bb Z}_n$ in circular way. We also denote the edge connecting $i$ to $i+1$ by $\xar{i}$. The left of Figure \ref{edge_pair} is the labels of a decagon. 

\begin{figure}[htp]
\centering
\begin{tikzpicture}[>=latex,scale=1]

\begin{scope}[xshift=-4cm]

\foreach \a in {0,...,9}
{
\draw
	(36*\a:1.2) -- (36*\a+36:1.2);
	
\node at (36*\a:1) {\footnotesize \a};
\node at (36*\a+18:1.35) {\footnotesize $\xar{\a}$};

}

\draw[->]
	(120:0.7) arc (120:420:0.7);

\end{scope}

\foreach \a in {0,...,9}
{
\foreach \b in {0,1}
\draw[xshift=4*\b cm]
	(36*\a:1.2) -- (36*\a+36:1.2);
	
\node at (36*\a:1) {\footnotesize \a};

}

\foreach \a/\b in {
	0/2, 1/4, 3/7, 5/8, 6/9}
\draw
	(36*\a+18:1.14) to[out=36*\a+198, in=36*\b+198] (36*\b+18:1.14);

\node at (0,-1.8) {\footnotesize $(\xar{0}\xar{2},\xar{1}\xar{4},\xar{3}\xar{7},\xar{5}\xar{8},\xar{6}\xar{9})$};

\begin{scope}[xshift=4cm]

\foreach \a/\b in {
	0/1, 2/-1, 1/1, 4/-1, 3/1, 7/-1, 5/1, 8/-1, 6/1, 9/-1}
\draw[->]
	(36*\a+18:1.14) -- ++(36*\a+108:0.1*\b);

\foreach \a/\b in {
	0/1, 2/1, 1/2, 4/2, 3/3, 7/3, 5/4, 8/4, 6/5, 9/5}
\node at (36*\a+18:1.35) {\footnotesize $a_{\b}$};	
	
\node at (0,-1.8) {\footnotesize $a_1a_2a_1^{-1}a_3a_2^{-1}a_4a_5a_3^{-1}a_4^{-1}a_5^{-1}$};

\end{scope}

\end{tikzpicture}
\caption{Label for polygon, and planar diagram.} 
\label{edge_pair}
\end{figure}

We denote a planar diagram for an orientable surface by 
\[
D=(\xar{i}_1\xar{j}_1,\xar{i}_2\xar{j}_2,\dots,\xar{i}_e\xar{j}_e),
\]
in which each edge appears exactly once. The orientability means that the edge pair $\xar{i}\xar{j}$ is identified in the opposing way, as the left of Figure \ref{edge_pair2}. 

For example, the middle of Figure \ref{edge_pair} is the planar diagram $D=(\xar{0}\xar{2},\xar{1}\xar{4},\xar{3}\xar{7},\xar{5}\xar{8},\xar{6}\xar{9})$, in which edge pairs are indicated by the chords connecting the middle points of the edges. The usual notation for the planar diagram labels the edge pairs as $a_1,a_2,a_3,a_4,a_5$, and the denotes the planar diagram as $a_1a_2a_1^{-1}a_3a_2^{-1}a_4a_5a_3^{-1}a_4^{-1}a_5^{-1}$.

\begin{figure}[htp]
\centering
\begin{tikzpicture}[>=latex,scale=1]

\foreach \a in {0,1}
{
\begin{scope}[xshift=4*\a cm]

\draw
	(-1,-0.5) -- (-1,0.5)
	(1,-0.5) -- (1,0.5);

\draw[dashed]
	(-1,0.5) to[out=90,in=90] (1,0.5)
	(-1,-0.5) to[out=-90,in=-90] (1,-0.5);

\node at (0.8,0) {\footnotesize $\xar{i}$};
\node at (-0.8,0) {\footnotesize $\xar{j}$};

\node at (0.8,-0.5) {\footnotesize $i$};
\node at (0.6,0.5) {\footnotesize $i\!+\!1$};
\node at (-0.8,0.5) {\footnotesize $j$};
\node at (-0.6,-0.5) {\footnotesize $j\!+\!1$};

\end{scope}
}


\draw[->]
	(1,0) -- ++(0,0.1);
\draw[->]
	(-1,0) -- ++(0,0.1);
	
\filldraw[fill=white]
	(1,-0.5) circle (0.08)
	(-1,-0.5) circle (0.08);

\fill 
	(1,0.5) circle (0.08)
	(-1,0.5) circle (0.08);

\node at (0,-1.5) {opposing};
\node at (0,0) {$+$};


\begin{scope}[xshift=4cm]

\draw[->]
	(1,0) -- ++(0,0.1);
\draw[->]
	(-1,0) -- ++(0,-0.1);
	
\fill 
	(1,-0.5) circle (0.08)
	(-1,0.5) circle (0.08);

\filldraw[fill=white]
	(1,0.5) circle (0.08)
	(-1,-0.5) circle (0.08);
	
\node at (0,-1.5) {twisted};
\node at (0,0) {$-$};

\end{scope}
		
\end{tikzpicture}
\caption{Opposing and twisted edge pairs.} 
\label{edge_pair2}
\end{figure}

Figure \ref{vertex_oppose} shows how an opposing pair of edges are glued together. We find that the pair induces adjacent corners $i+1$ and $j$ at a vertex $\bullet$ of the tiling, and another adjacent corners $j+1$ and $i$ at another vertex $\circ$. We denote the process by the following, where the pairs on the left are ordered
\begin{equation}\label{pair2vertex_o}
\xar{i}\xar{j}
\implies \bullet = (i+1)j\cdots,\quad
\xar{j}\xar{i}
\implies \circ = (j+1)i\cdots.
\end{equation}

\begin{figure}[htp]
\centering
\begin{tikzpicture}[>=latex,scale=1]

\draw
	(0,-1.2) -- (0,1.2)
	(-1.8,-1.2) -- (1.8,-1.2)
	(-1.8,1.2) -- (1.8,1.2);

\fill 
	(0,-1.2) circle (0.08);
\filldraw[fill=white]
	(0,1.2) circle (0.08);
		
\draw[yshift=1.2 cm, <-]
	(15:0.7) arc (15:-195:0.7);
	
\draw[yshift=-1.2 cm, ->]
	(-15:0.7) arc (-15:195:0.7);

\foreach \a in {1,-1}
\draw[xshift=1.4*\a cm, ->]
	(120:0.3) arc (120:420:0.3);
		
\node at (-0.2,0) {\footnotesize $\xar{j}$};
\node at (0.2,0) {\footnotesize $\xar{i}$};

\node at (-1.2,-0.95) {\footnotesize $\overline{j\!-\!1}$};
\node at (-1.2,0.95) {\footnotesize $\overline{j\!+\!1}$};

\node at (-0.3,1) {\scriptsize $j\!+\!1$};
\node at (-0.2,-1) {\scriptsize $j$};

\node at (1.2,-0.95) {\footnotesize $\overline{i\!+\!1}$};
\node at (1.2,0.95) {\footnotesize $\overline{i\!-\!1}$};

\node at (0.33,-1) {\scriptsize $i\!+\!1$};
\node at (0.2,1) {\scriptsize $i$};
		
\end{tikzpicture}
\caption{Opposing edge pairs inducing adjacent corners at vertices.} 
\label{vertex_oppose}
\end{figure}

We remark that, although the edge pair $\xar{i}\xar{j}$ in a planar diagram has no order among $\xar{i}$ and $\xar{j}$, we do care about the order when they induce adjacent corners at vertices in \eqref{pair2vertex_o}. The two formulas are actually the same, after exchanging $i$ and $j$.

We may further construct the vertex $\bullet$ by considering which edge is paired with $\xar{j-1}$. This adds a corner adjacent to (in the counterclockwise direction) the corner $j$ at $\bullet$. Continuing the process eventually comes back to the corner $i+1$, and completes the vertex $\bullet$. The same can be carried out for the vertex $\circ$. For example, we rewrite the formula \eqref{pair2vertex_o} as 
\[
i\cdots
\overset{\scriptsize \xar{i-1}\xar{j}}{=\joinrel=\joinrel=} 
ij\cdots,
\]
and calculate vertices from the decagon planar diagram $D=(\xar{0}\xar{2},\xar{1}\xar{4},\xar{3}\xar{7},\xar{5}\xar{8},\xar{6}\xar{9})$ in Figure \ref{edge_pair}
\begin{align*}
\bullet
=0\cdots
&\overset{\scriptsize \xar{9}\xar{6}}{=\joinrel=} 
06\cdots
\overset{\scriptsize \xar{5}\xar{8}}{=\joinrel=} 
068\cdots
\overset{\scriptsize \xar{7}\xar{3}}{=\joinrel=} 
0683\cdots
\overset{\scriptsize \xar{2}\xar{0}}{=\joinrel=} 
0683, \\
\circ
=1\cdots
&\overset{\scriptsize \xar{0}\xar{2}}{=\joinrel=} 
12\cdots
\overset{\scriptsize \xar{1}\xar{4}}{=\joinrel=} 
124\cdots
\overset{\scriptsize \xar{3}\xar{7}}{=\joinrel=} 
1247\cdots
\overset{\scriptsize \xar{6}\xar{9}}{=\joinrel=} 
12479\cdots \\
&\overset{\scriptsize \xar{8}\xar{5}}{=\joinrel=} 
124795\cdots 
\overset{\scriptsize \xar{4}\xar{1}}{=\joinrel=} 
124795.
\end{align*}
The calculation shows that the single tile tiling has two vertices $\bullet=0683$ and $\circ=124795$ in the middle of Figure \ref{pd_vs}.

\begin{figure}[htp]
\centering
\begin{tikzpicture}[>=latex,scale=1]

\draw[gray, very thick]
	(0:0.9) to[out=180+36*0, in=180+36*6] (36*6:0.9) to[out=180+36*6, in=180+36*8] (36*8:0.9) to[out=180+36*8, in=180+36*3] (36*3:0.9) to[out=180+36*3, in=180+36*0] (36*0:0.9);
	
\draw[gray, very thick, ->]
	(0,0) -- ++(108:0.3);

\foreach \a in {0,...,9}
{
\draw
	(36*\a:1.2) -- (36*\a+36:1.2);
	
\node at (36*\a:1) {\footnotesize \a};

}

\foreach \a in {0,3,6,8}
\fill 
	(36*\a:1.2) circle (0.05);

\foreach \a in {1,2,4,5,7,9}
\filldraw[fill=white]
	(36*\a:1.2) circle (0.05);

\foreach \a/\b in {
	0/2, 1/4, 3/7, 5/8, 6/9}
\draw
	(36*\a+18:1.14) to[out=36*\a+198, in=36*\b+198] (36*\b+18:1.14);
	
\draw[->, very thick]
	(1.8,0) -- ++(1,0);

\begin{scope}[xshift=4cm]

\fill 
	(0,0) circle (0.08); 

\draw[<-] (240:0.3) arc (240:-60:0.3);

\foreach \a/\b in {0/0, 1/3, 2/8, 3/6}
\node at (-90*\a:0.5) {\footnotesize \b};

\node at (0,-1) {\footnotesize 0683};

\end{scope}

\begin{scope}[xshift=6cm]

\filldraw[fill=white]
	(0,0) circle (0.08); 

\draw[<-] (240:0.3) arc (240:-60:0.3);

\foreach \a/\b in {0/1, 1/5, 2/9, 3/7, 4/4, 5/2}
\node at (-60*\a:0.5) {\small \b};

\node at (0,-1) {\small 124795};

\end{scope}

\begin{scope}[xshift=10cm]

\draw
	(50:1.2) -- (80:1.2)
	(150:1.2) -- (180:1.2)
	(65:1.16) to[out=245, in=-15] (165:1.16);

\draw[dashed]
	(80:1.2) arc (80:150:1.2)
	(180:1.2) arc (180:410:1.2);

\filldraw[fill=white]
	(50:1.2) circle (0.05)
	(180:1.2) circle (0.05);
	
\fill
	(80:1.2) circle (0.05)
	(150:1.2) circle (0.05);

\draw[gray, very thick, <-]
	(50:1.1) to[out=180+50, in=0] (180:1.1);

\draw[gray, very thick, ->]
	(80:1.1) to[out=180+80, in=180+150] (150:1.1);
	
\end{scope}
	
\end{tikzpicture}
\caption{Planar diagram to vertices.} 
\label{pd_vs}
\end{figure}

We can easily read the vertices from the chords in the picture of a planar diagram. On the right of Figure \ref{pd_vs}, we see the two sides of a chord connect pairs of vertices by thick gray lines. If we connect vertices of the planar diagram by following this rule, we get directed closed curves. Then each closed curve is a vertex. For example, the thick gray quadrilateral on the left indicates the vertex 0683.

In general, a planar diagram $D$ divides all the corners into a disjoint union of circularly ordered subsets $V_1,V_2,\dots,V_v$. The vertex set in the example above consists of $V_1=0683$ and $V_2=124795$. By circularly ordered, we mean rotations represent the same vertex
\[
0683=6830=8306=3068.
\]
Each subset is a {\em vertex} of the tiling ${\mc T}_D$, and the collection of the circularly ordered subsets is a {\em vertex set}. 

The degree of a vertex is the number of corners in the subset. In the example above, the degrees of $\bullet$ and $\circ$ are respectively 4 and 6. If $v$ is the number of subsets, then the Euler number of the surface is $\chi=v-e+1$. The example has $\chi=2-5+1=-2$. This means the surface $S_D=2{\bb T}^2$.

Conversely, we can recover the planar diagram from a vertex set. Specifically, we reverse the formula \eqref{pair2vertex_o} to get the formula from vertices to edge pairs
\begin{equation}\label{vertex2pair_o}
ij\cdots \implies \xar{i-1}\xar{j}.
\end{equation}
For example, the vertex $0683$ can be expressed as $06\cdots,68\cdots,83\cdots,30\cdots$, which give the edge pairs
\[
06\cdots \implies
\xar{9}\xar{6},\quad
68\cdots \implies 
\xar{5}\xar{8},\quad
83\cdots \implies 
\xar{7}\xar{3},\quad
30\cdots \implies 
\xar{2}\xar{0}.
\]
By exchanging the order in the edge pairs, we find that the following should also be vertices
\[
\xar{0}\xar{2}
\implies 12\cdots,\quad
\xar{3}\xar{7}
\implies 47\cdots,\quad
\xar{8}\xar{5}
\implies 95\cdots,\quad
\xar{6}\xar{9}
\implies 79\cdots.
\]
Then we find $12\cdots$ and $4795\cdots$ should be vertices. Moreover, the four known edge pairs imply the fifth edge pair $\xar{1}\xar{4}$, which further implies $24\cdots$ and $51\cdots$ are vertices. Therefore $124795$ is a vertex. 

We conclude that the vertex $0683$ implies there is only one more vertex $124795$, and also uniquely determines the planar diagram $(\xar{0}\xar{2},\xar{1}\xar{4},\xar{3}\xar{7},\xar{5}\xar{8},\xar{6}\xar{9})$.

We see that planar diagrams are equivalent to vertex sets. In a vertex set, we require all the circularly ordered subsets to have at least three elements, and the edge pairs induced by all the circularly ordered subsets are compatible. The requirement is the same as the following.

\begin{proposition}\label{avs}
A disjoint union of circularly ordered subsets $V_1,V_2,\dots,V_v$ of ${\bb Z}_n$ is a vertex set (i.e., induced from a planar diagram) if and only if the following are satisfied:
\begin{itemize}
\item Each (circularly ordered) subset has at least three corners.
\item $i(i-1)\cdots$ is not a subset.
\item If $ij\cdots$ is a subset, then $(j+1)(i-1)\cdots$ is a subset.
\end{itemize}
\end{proposition}

The subsets will be vertices. However, we still call them subsets instead of vertices in the statement of the lemma and the if part of the proof, to avoid suggesting that they are already vertices in some tiling.

We remark that, in addition to the second condition in Lemma \ref{avs}, another useful condition is that $i(i-2)\cdots$ is not a subset. The reason is that applying the third condition to $i(i-2)\cdots$ gives a subset $(i-1)(i-1)\cdots$, a contradiction.  

\begin{proof}
Suppose we have a vertex set. Then the first condition is satisfied because each vertex has degree $\ge 3$. The second condition is satisfied because adjacent corners $i(i-1)$ at a vertex would imply $\xar{i-1}\xar{i-1}$ is an edge pair, a contradiction. The third condition is satisfied because a vertex $ij\cdots$ implies an edge pair $\xar{i-1}\xar{j}$, and then the edge pair $\xar{j}\xar{i-1}$ implies the vertex $(j+1)(i-1)\cdots$. 

Conversely, suppose a disjoint union satisfies three conditions. We consider the map $f(ij)=\xar{i-1}\xar{j}$ in \eqref{vertex2pair_o}, from the set $X$ of adjacent corners at vertices to the set $Y$ of induced edge pairs. Note that the pairs in $X$ are ordered, and the pair in $Y$ are not ordered.

The third condition means that, if $f((i+1)j)=\xar{i}\xar{j}$ for some $(i+1)j\in X$, then we also have $f((j+1)i)=\xar{j}\xar{i}$ for another $(j+1)i\in X$. By the second condition, the adjacent corners $(i+1)j$ and $(j+1)i$ at vertices are different pairs. Therefore the map $f$ is at least two to one everywhere.

On the other hand, since every corner $j\in {\bb Z}_n$ appears in a subset $ij\cdots$, the corresponding edge $\xar{j}$ appears in an induced edge pair $f(ij)=\xar{i-1}\xar{j}$. Therefore there are at least $\frac{n}{2}=e$ edge pairs in $Y$. Since $X$ has $n=2e$ elements, and $f$ is at least two to one everywhere, we conclude $f$ is exactly two to one everywhere. This implies $Y$ has exactly $e$ elements, and is a planar diagram.
\end{proof}

Two planar diagrams represent the same tiling if and only if they are related by relabelling the $n$-gon. There are two types of relabelling ($c\in {\bb Z}_n$ is a constant)
\begin{itemize}
\item Rotation: $i\mapsto c+i$ for corners, and $\bar{i}\mapsto \overline{c+i}$ for edges.
\item Reversion: $i\mapsto c-i$ for corners, and $\bar{i}\mapsto \overline{c-i-1}$ for edges. After applying the relabelling to corners, the circular order of vertices needs to be reversed.
\end{itemize}
The following is an example of the reversion of the planar diagram for $c=4$
\[
(\xar{0}\xar{2},\xar{1}\xar{4},\xar{3}\xar{7},\xar{5}\xar{8},\xar{6}\xar{9})
\xrightarrow{\overline{4-i-1}}
(\xar{3}\xar{1},\xar{2}\xar{9},\xar{0}\xar{6},\xar{8}\xar{5},\xar{7}\xar{4})
=(\xar{0}\xar{6},\xar{1}\xar{3},\xar{2}\xar{9},\xar{4}\xar{7},\xar{5}\xar{8}).
\]
This is equivalent to the following reversion of the vertices
\[
(0683, 124795) 
\xrightarrow{4-i} 
(4861, 320759) 
\xrightarrow{\text{reverse}}
(1684, 957023)
=(023957,1684).
\]

We may create computer program to find all the single tile tilings of a fixed orientable surface of Euler characteristic $<0$, by congruent $n$-gons with $n\ge 7$. The program can use either the code for planar diagram or the code for vertex set. 

The constraints on the code for planar diagram are the following
\begin{itemize}
\item No degree one vertex: $\xar{i}\xar{i+1}$ is not an edge pair.
\item No degree two vertex: If $\xar{i}\xar{j}$ is an edge pair, then $\xar{i+1}\xar{j-1}$ is not an edge pair.
\end{itemize} 
The constraints are rather simple. The complication is the determination of the surface. This can only be done by calculating all the vertices, and then using the number of vertices to calculate the Euler number. If we wish to find all the single tile tilings of an orientable surface, then we need to find all the planar diagrams that yield specific number of vertices.

The advantage of using the code for vertices is the explicit number of vertices, which corresponds to specific surface. The complication lies in the constraints in Proposition \ref{avs}.  

In a tiling of $2{\bb T}^2$ by a single $n$-gon, we know $n=2e$ with $4\le e\le 9$. The number of vertices in the tiling is $v=e-3$. For $3{\bb T}^2$, we know $n=2e$ with $6\le e\le 15$, and $v=e-5$. For $4{\bb T}^2$, we know $n=2e$ with $8\le e\le 21$, and $v=e-7$. Table \ref{T2tiling} shows the numbers of tilings. Figures \ref{2T2_n8} and \ref{2T2_n16} show all the 122 tilings of $2{\bb T}^2$ by a single $n$-gon with $n\ge 7$. We can easily get the equivalent vertex sets by following the rule in Figure \ref{pd_vs}.

\begin{table}[htp]
	\centering
	\begin{tabular}{|c|c||c|c|c|c|c|c|c|}
	\hline
	\multirow{2}{*}{$2{\bb T}^2$}
	& polygon size & 8 & 10 & 12 & 14 & 16 & 18 \\
	\cline{2-8}
	& number       & 4 & 18 & 34 & 38 & 20 & 8 \\
	\hline
	\hline
	\multirow{4}{*}{$3{\bb T}^2$}
	& polygon size &&& 12 & 14 & 16 & 18 \\
	\cline{2-8}
	& number       &&& 82 & 1022 & 5741 & 18281 \\
	\cline{2-8} 
	& polygon size & 20 & 22 & 24 & 26 & 28 & 30 \\
	\cline{2-8}
	& number       & 36232 & 46784 & 32296 & 20978 & 6396 & 927 \\
	\hline \hline
	\multirow{2}{*}{$4{\bb T}^2$}
	& polygon size & 16 & 18 & 20 & 22 & $\cdots$ & 42 \\
	\cline{2-8}
	& number       & 7258 & 175136 & 1785661 & 10404513 & $\cdots$ & 676445 \\
	\hline
	\end{tabular}
\caption{Number of single tile tilings of orientable hyperbolic surfaces of small genus. }
\label{T2tiling}
\end{table}

\begin{figure}[htp]
\centering
\begin{tikzpicture}[>=latex,scale=1]


\begin{scope}[shift={(4.5cm, 5.6cm)}]

\foreach \x in {0,...,3}
\foreach \a in {0,...,7}
{
\begin{scope}[xshift=1.8*\x cm]

\draw
	(45*\a-22.5:0.8) -- (45*\a+22.5:0.8);

\fill[red]
	(45*\a+22.5:0.8) circle (0.05);
	
\end{scope}
}

\foreach \x in {1,...,4}
\node[gray] at (-1.8+1.8*\x, 0) { 8.\x};

\foreach \a/\b/\x in {
	0/4/0, 1/5/0, 2/6/0, 3/7/0,
	0/2/1, 3/6/1, 4/7/1, 1/5/1,
	0/2/2, 1/3/2, 4/6/2, 5/7/2,
	0/2/3, 5/7/3, 1/4/3, 3/6/3}
\draw[xshift=1.8*\x cm]
	(45*\a+45:0.74) to[out=45*\a+225, in=45*\b+225] (45*\b+45:0.74);
	
\end{scope}


\begin{scope}[yshift=3.8cm]

\foreach \x in {0,...,8}
\foreach \y in {0,1}
\foreach \a in {0,...,9}
\draw[shift={(1.8*\x cm,-1.8*\y cm)}]
	(36*\a:0.8) -- (36*\a+36:0.8);

\begin{scope}[gray!70]

\foreach \x in {1,...,9}
\node at (-1.8+1.8*\x, 0) {10.\x};

\foreach \x in {10,...,18}
\node at (-18+1.8*\x, -1.8) {10.\x};

\end{scope}

\foreach \a/\b/\x/\y in { 
	0/5/0/0, 1/6/0/0, 2/7/0/0, 3/8/0/0, 4/9/0/0,  
	6/8/1/0, 2/7/1/0, 3/9/1/0, 0/4/1/0, 1/5/1/0, 
	1/3/2/0, 2/7/2/0, 4/9/2/0, 0/5/2/0, 6/8/2/0, 
	0/3/3/0, 1/4/3/0, 2/7/3/0, 5/8/3/0, 6/9/3/0, 
	0/8/4/0, 3/9/4/0, 1/5/4/0, 2/7/4/0, 4/6/4/0, 
	1/9/5/0, 0/4/5/0, 2/7/5/0, 3/5/5/0, 6/8/5/0, 
	0/8/6/0, 1/9/6/0, 2/7/6/0, 3/5/6/0, 4/6/6/0, 
	1/9/7/0, 0/3/7/0, 2/7/7/0, 4/6/7/0, 5/8/7/0, 
	0/2/8/0, 1/4/8/0, 3/7/8/0, 5/8/8/0, 6/9/8/0, 
	0/3/0/1, 1/8/0/1, 2/5/0/1, 4/7/0/1, 6/9/0/1, 
	1/8/1/1, 3/9/1/1, 0/6/1/1, 2/5/1/1, 4/7/1/1, 
	1/8/2/1, 4/9/2/1, 0/5/2/1, 2/6/2/1, 3/7/2/1, 
	0/2/3/1, 1/5/3/1, 3/6/3/1, 4/8/3/1, 7/9/3/1,
	0/2/4/1, 1/4/4/1, 3/6/4/1, 5/8/4/1, 7/9/4/1, 
	0/8/5/1, 1/9/5/1, 2/5/5/1, 3/6/5/1, 4/7/5/1, 
	0/8/6/1, 1/9/6/1, 2/4/6/1, 3/6/6/1, 5/7/6/1, 
	0/3/7/1, 1/5/7/1, 2/7/7/1, 4/8/7/1, 6/9/7/1, 
	0/2/8/1, 1/5/8/1, 3/7/8/1, 4/8/8/1, 6/9/8/1
	}
\draw[xshift=1.8*\x cm, yshift=-1.8*\y cm]
	(36*\a+18:0.76) to[out=36*\a+198, in=36*\b+198] (36*\b+18:0.76);

\foreach \a/\x/\y in {
	0/0/0, 2/0/0, 4/0/0, 6/0/0, 8/0/0,
	0/1/0, 2/1/0, 3/1/0, 5/1/0, 7/1/0, 8/1/0,
	2/2/0, 3/2/0, 7/2/0, 8/2/0,
	0/3/0, 2/3/0, 4/3/0, 6/3/0, 8/3/0,
	1/4/0, 2/4/0, 5/4/0, 6/4/0, 8/4/0,
	2/5/0, 3/5/0, 6/5/0, 7/5/0, 8/5/0, 9/5/0,
	0/6/0, 1/6/0, 2/6/0, 8/6/0, 9/6/0,
	2/7/0, 5/7/0, 6/7/0, 8/7/0, 9/7/0,
	1/8/0, 2/8/0, 4/8/0, 5/8/0, 7/8/0, 9/8/0,
	0/0/1, 2/0/1, 4/0/1, 6/0/1, 8/0/1,
	1/1/1, 2/1/1, 4/1/1, 6/1/1, 8/1/1, 9/1/1,
	0/2/1, 2/2/1, 3/2/1, 4/2/1, 6/2/1, 7/2/1, 8/2/1,
	1/3/1, 2/3/1, 4/3/1, 5/3/1, 6/3/1, 8/3/1, 9/3/1,
	1/4/1, 2/4/1, 4/4/1, 5/4/1, 6/4/1, 8/4/1, 9/4/1,
	0/5/1, 1/5/1, 2/5/1, 4/5/1, 6/5/1, 8/5/1, 9/5/1,
	0/6/1, 1/6/1, 2/6/1, 5/6/1, 8/6/1, 9/6/1,
	2/7/1, 5/7/1, 8/7/1,
	0/8/1, 1/8/1, 2/8/1, 3/8/1, 5/8/1, 6/8/1, 8/8/1
	}	
\fill[red, shift={(1.8*\x cm,-1.8*\y cm)}]
	(36*\a:0.8) circle (0.05);

\foreach \a/\x/\y in {
	1/0/0, 3/0/0, 5/0/0, 7/0/0, 9/0/0,
	1/1/0, 4/1/0, 6/1/0, 9/1/0,
	0/2/0, 1/2/0, 4/2/0, 5/2/0, 6/2/0, 9/2/0,
	1/3/0, 3/3/0, 5/3/0, 7/3/0, 9/3/0,
	0/4/0, 3/4/0, 4/4/0, 7/4/0, 9/4/0,
	0/5/0, 1/5/0, 4/5/0, 5/5/0,
	3/6/0, 4/6/0, 5/6/0, 6/6/0, 7/6/0,
	0/7/0, 1/7/0, 3/7/0, 4/7/0, 7/7/0,
	0/8/0, 3/8/0, 6/8/0, 8/8/0,
	1/0/1, 3/0/1, 5/0/1, 7/0/1, 9/0/1,
	0/1/1, 3/1/1, 5/1/1, 7/1/1,
	1/2/1, 5/2/1, 9/2/1,
	0/3/1, 3/3/1, 7/3/1,
	0/4/1, 3/4/1, 7/4/1,
	3/5/1, 5/5/1, 7/5/1,
	3/6/1, 4/6/1, 6/6/1, 7/6/1,
	0/7/1, 1/7/1, 3/7/1, 4/7/1, 6/7/1, 7/7/1, 9/7/1,
	4/8/1, 7/8/1, 9/8/1
	}	
\fill[blue, shift={(1.8*\x cm,-1.8*\y cm)}]
	(36*\a:0.8) circle (0.05);
	
\end{scope}


\foreach \x in {0,...,8}
\foreach \y in {0,...,2}
\foreach \a in {0,...,11}
\draw[shift={(1.8*\x cm,-1.8*\y cm)}]
	(30*\a-15:0.8) -- (30*\a+15:0.8);

\foreach \x in {1,...,7}
\foreach \a in {0,...,11}
\draw[shift={(1.8*\x cm,-1.8*3 cm)}]
	(30*\a-15:0.8) -- (30*\a+15:0.8);
	
\begin{scope}[gray!70]

\foreach \x in {1,...,9}
\node at (-1.8+1.8*\x, 0) {12.\x};

\foreach \x in {10,...,18}
\node at (-18+1.8*\x, -1.8) {12.\x};

\foreach \x in {19,...,27}
\node at (-34.2+1.8*\x, -3.6) {12.\x};

\foreach \x in {28,...,34}
\node at (-48.6+1.8*\x, -5.4) {12.\x};

\end{scope}

\foreach \a/\b/\x/\y in {
0/3/0/0, 1/10/0/0, 2/5/0/0, 4/7/0/0, 6/9/0/0, 8/11/0/0, 
2/11/1/0, 0/5/1/0, 1/8/1/0, 3/6/1/0, 4/9/1/0, 7/10/1/0,
2/11/2/0, 0/6/2/0, 1/7/2/0, 3/9/2/0, 4/10/2/0, 5/8/2/0, 
5/8/3/0, 6/10/3/0, 3/7/3/0, 1/9/3/0, 2/11/3/0, 0/4/3/0,
2/11/4/0, 0/4/4/0, 1/9/4/0, 3/6/4/0, 5/8/4/0, 7/10/4/0, 
2/11/5/0, 0/5/5/0, 1/8/5/0, 3/7/5/0, 4/9/5/0, 6/10/5/0, 
1/10/6/0, 6/11/6/0, 0/3/6/0, 2/7/6/0, 4/8/6/0, 5/9/6/0,  
1/4/7/0, 2/7/7/0, 3/8/7/0, 5/10/7/0, 6/11/7/0, 0/9/7/0,  
1/4/8/0, 2/6/8/0, 3/10/8/0, 5/8/8/0, 7/11/8/0, 0/9/8/0,
0/10/0/1, 2/11/0/1, 1/3/0/1, 4/6/0/1, 5/8/0/1, 7/9/0/1, 
1/10/1/1, 2/11/1/1, 0/3/1/1, 4/7/1/1, 5/8/1/1, 6/9/1/1,
1/11/2/1, 0/2/2/1, 3/6/2/1, 4/8/2/1, 5/9/2/1, 7/10/2/1, 
0/10/3/1, 2/11/3/1, 1/3/3/1, 4/7/3/1, 5/8/3/1, 6/9/3/1,  
1/11/4/1, 0/2/4/1, 3/10/4/1, 4/7/4/1, 5/8/4/1, 6/9/4/1,
1/11/5/1, 0/2/5/1, 3/10/5/1, 4/6/5/1, 5/8/5/1, 7/9/5/1,
5/8/6/1, 2/6/6/1, 7/11/6/1, 1/9/6/1, 3/10/6/1, 0/4/6/1,  
9/11/7/1, 3/10/7/1, 0/6/7/1, 1/7/7/1, 2/4/7/1, 5/8/7/1,  
2/11/8/1, 0/6/8/1, 1/7/8/1, 3/8/8/1, 4/9/8/1, 5/10/8/1, 
1/11/0/2, 0/3/0/2, 2/8/0/2, 4/7/0/2, 5/9/0/2, 6/10/0/2,  
0/10/1/2, 3/11/1/2, 1/5/1/2, 2/8/1/2, 4/7/1/2, 6/9/1/2,
4/6/2/2, 5/10/2/2, 0/7/2/2, 2/8/2/2, 9/11/2/2, 1/3/2/2, 
3/5/3/2, 4/9/3/2, 6/11/3/2, 0/7/3/2, 2/8/3/2, 1/10/3/2, 
0/3/4/2, 1/5/4/2, 2/8/4/2, 4/9/4/2, 6/10/4/2, 7/11/4/2,
0/10/5/2, 3/11/5/2, 1/6/5/2, 2/8/5/2, 4/7/5/2, 5/9/5/2,  
0/2/6/2, 1/4/6/2, 3/7/6/2, 5/9/6/2, 6/10/6/2, 8/11/6/2, 
0/2/7/2, 1/4/7/2, 3/8/7/2, 5/7/7/2, 6/10/7/2, 9/11/7/2, 
0/2/8/2, 1/5/8/2, 3/7/8/2, 4/9/8/2, 6/10/8/2, 8/11/8/2, 
1/3/1/3, 2/8/1/3, 4/6/1/3, 5/11/1/3, 7/9/1/3, 0/10/1/3, 
0/10/2/3, 5/11/2/3, 1/7/2/3, 2/8/2/3, 3/9/2/3, 4/6/2/3, 
0/4/3/3, 1/6/3/3, 2/8/3/3, 3/10/3/3, 5/9/3/3, 7/11/3/3, 
1/11/4/3, 0/4/4/3, 2/8/4/3, 3/5/4/3, 6/9/4/3, 7/10/4/3, 
1/3/5/3, 2/8/5/3, 4/7/5/3, 5/10/5/3, 6/11/5/3, 0/9/5/3, 
0/2/6/3, 1/6/6/3, 3/7/6/3, 4/9/6/3, 5/10/6/3, 8/11/6/3, 
0/3/7/3, 1/4/7/3, 2/7/7/3, 5/9/7/3, 6/10/7/3, 8/11/7/3
}
\draw[xshift=1.8*\x cm, yshift=-1.8*\y cm]
	(30*\a+30:0.77) to[out=30*\a+210, in=30*\b+210] (30*\b+30:0.77);

\foreach \a/\x/\y in {
	0/0/0, 2/0/0, 4/0/0, 6/0/0, 8/0/0, 10/0/0,
	0/1/0, 4/1/0, 8/1/0,
	0/2/0, 5/2/0, 7/2/0,
	0/3/0, 3/3/0, 6/3/0, 9/3/0,
	0/4/0, 3/4/0, 5/4/0, 7/4/0, 9/4/0,
	0/5/0, 4/5/0, 8/5/0,
	0/6/0, 2/6/0, 6/6/0, 10/6/0,
	0/7/0, 2/7/0, 4/7/0, 6/7/0, 8/7/0, 10/7/0,
	0/8/0, 4/8/0, 8/8/0,
	0/0/1, 3/0/1, 6/0/1, 9/0/1,
	0/1/1, 2/1/1, 10/1/1,
	0/2/1, 1/2/1, 2/2/1, 6/2/1, 10/2/1, 11/2/1,
	0/3/1, 3/3/1, 5/3/1, 7/3/1, 9/3/1,
	0/4/1, 1/4/1, 2/4/1, 10/4/1, 11/4/1,
	0/5/1, 1/5/1, 2/5/1, 10/5/1, 11/5/1,
	0/6/1, 3/6/1, 9/6/1,
	0/7/1, 5/7/1, 7/7/1,
	0/8/1, 3/8/1, 5/8/1, 7/8/1, 9/8/1,
	1/0/2, 5/0/2, 8/0/2, 10/0/2,
	1/1/2, 4/1/2, 6/1/2, 8/1/2,
	1/2/2, 2/2/2, 7/2/2, 8/2/2, 11/2/2,
	1/3/2, 3/3/2, 4/3/2, 8/3/2, 9/3/2,
	1/4/2, 4/4/2, 8/4/2,
	1/5/2, 5/5/2, 8/5/2,
	0/6/2, 1/6/2, 3/6/2, 4/6/2, 6/6/2, 9/6/2,
	0/7/2, 1/7/2, 3/7/2, 4/7/2, 7/7/2,
	0/8/2, 1/8/2, 4/8/2, 5/8/2, 8/8/2, 10/8/2,
	1/1/3, 2/1/3, 7/1/3, 8/1/3,
	1/2/3, 3/2/3, 6/2/3, 8/2/3,
	1/3/3, 5/3/3, 8/3/3,
	1/4/3, 6/4/3, 8/4/3, 10/4/3,
	0/5/3, 1/5/3, 2/5/3, 3/5/3, 7/5/3, 8/5/3,
	0/6/3, 1/6/3, 3/6/3, 5/6/3, 6/6/3, 9/6/3,
	1/7/3, 3/7/3, 7/7/3, 11/7/3
	}	
\fill[red, shift={(1.8*\x cm,-1.8*\y cm)}]
	(30*\a+45:0.8) circle (0.05);

\foreach \a/\x/\y in {
	1/0/0, 5/0/0, 9/0/0,
	1/1/0, 3/1/0, 5/1/0, 7/1/0, 9/1/0, 11/1/0,
	1/2/0, 6/2/0, 11/2/0,
	1/3/0, 4/3/0, 8/3/0, 11/3/0,
	1/4/0, 4/4/0, 8/4/0, 11/4/0,
	1/5/0, 2/5/0, 5/5/0, 7/5/0, 10/5/0, 11/5/0,
	1/6/0, 4/6/0, 7/6/0, 9/6/0,
	1/7/0, 3/7/0, 7/7/0,
	1/8/0, 3/8/0, 6/8/0, 9/8/0, 11/8/0,
	1/0/1, 2/0/1, 10/0/1, 11/0/1,
	4/1/1, 6/1/1, 8/1/1,
	3/2/1, 5/2/1, 8/2/1,
	1/3/1, 2/3/1, 10/3/1, 11/3/1,
	3/4/1, 5/4/1, 7/4/1, 9/4/1,
	3/5/1, 6/5/1, 9/5/1,
	1/6/1, 4/6/1, 6/6/1, 8/6/1, 11/6/1,
	1/7/1, 4/7/1, 6/7/1, 8/7/1, 11/7/1,
	1/8/1, 6/8/1, 11/8/1,
	0/0/2, 2/0/2, 3/0/2, 7/0/2, 11/0/2,
	2/1/2, 3/1/2, 7/1/2, 10/1/2, 11/1/2,
	0/2/2, 3/2/2, 6/2/2,
	2/3/2, 5/3/2, 7/3/2, 11/3/2,
	0/4/2, 2/4/2, 5/4/2, 7/4/2, 10/4/2,
	2/5/2, 3/5/2, 7/5/2, 10/5/2, 11/5/2,
	2/6/2, 7/6/2, 11/6/2,
	2/7/2, 8/7/2, 11/7/2,
	2/8/2, 7/8/2, 11/8/2,
	4/1/3, 5/1/3, 10/1/3, 11/1/3,
	0/2/3, 2/2/3, 7/2/3, 9/2/3,
	2/3/3, 7/3/3, 10/3/3,
	2/4/3, 5/4/3, 7/4/3, 9/4/3,
	4/5/3, 6/5/3, 10/5/3,
	2/6/3, 7/6/3, 11/6/3,
	0/7/3, 2/7/3, 4/7/3, 6/7/3, 9/7/3
	}	
\fill[blue, shift={(1.8*\x cm,-1.8*\y cm)}]
	(30*\a+45:0.8) circle (0.05);

\foreach \a/\x/\y in {
	3/0/0, 7/0/0, 11/0/0,
	2/1/0, 6/1/0, 10/1/0,
	2/2/0, 3/2/0, 4/2/0, 8/2/0, 9/2/0, 10/2/0,
	2/3/0, 5/3/0, 7/3/0, 10/3/0,
	2/4/0, 6/4/0, 10/4/0,
	3/5/0, 6/5/0, 9/5/0,
	3/6/0, 5/6/0, 8/6/0, 11/6/0,
	5/7/0, 9/7/0, 11/7/0,
	2/8/0, 5/8/0, 7/8/0, 10/8/0,
	4/0/1, 5/0/1, 7/0/1, 8/0/1,
	1/1/1, 3/1/1, 5/1/1, 7/1/1, 9/1/1, 11/1/1,
	4/2/1, 7/2/1, 9/2/1,
	4/3/1, 6/3/1, 8/3/1,
	4/4/1, 6/4/1, 8/4/1,
	4/5/1, 5/5/1, 7/5/1, 8/5/1,
	2/6/1, 5/6/1, 7/6/1, 10/6/1,
	2/7/1, 3/7/1, 9/7/1, 10/7/1,
	2/8/1, 4/8/1, 8/8/1, 10/8/1,
	4/0/2, 6/0/2, 9/0/2,
	0/1/2, 5/1/2, 9/1/2,
	4/2/2, 5/2/2, 9/2/2, 10/2/2,
	0/3/2, 6/3/2, 10/3/2,
	3/4/2, 6/4/2, 9/4/2, 11/4/2,
	0/5/2, 4/5/2, 6/5/2, 9/5/2,
	5/6/2, 8/6/2, 10/6/2,
	5/7/2, 6/7/2, 9/7/2, 10/7/2,
	3/8/2, 6/8/2, 9/8/2,
	0/1/3, 3/1/3, 6/1/3, 9/1/3,
	4/2/3, 5/2/3, 10/2/3, 11/2/3,
	0/3/3, 3/3/3, 4/3/3, 6/3/3, 9/3/3, 11/3/3,
	0/4/3, 3/4/3, 4/4/3, 11/4/3,
	5/5/3, 9/5/3, 11/5/3,
	4/6/3, 8/6/3, 10/6/3,
	5/7/3, 8/7/3, 10/7/3
	}	
\fill[green, shift={(1.8*\x cm,-1.8*\y cm)}]
	(30*\a+45:0.8) circle (0.05);

\begin{scope}[shift={(0.2cm, -7.4cm)}]


\foreach \a in {0,...,13}
{
\foreach \x in {0,...,7}
\foreach \y in {0,...,3}
\draw[shift={(2*\x cm,-2*\y cm)}]
	(25.714*\a:0.9) -- (25.714*\a+25.714:0.9);	

\foreach \x in {1,...,6}
\draw[shift={(2*\x cm,-2*4 cm)}]
	(25.714*\a:0.9) -- (25.714*\a+25.714:0.9);	
}

\begin{scope}[gray!70]

\foreach \x in {1,...,8}
\node at (-2+2*\x, 0) { 14.\x};

\foreach \x in {9,...,16}
\node at (-18+2*\x, -2) { 14.\x};

\foreach \x in {17,...,24}
\node at (-34+2*\x, -4) { 14.\x};

\foreach \x in {25,...,32}
\node at (-50+2*\x, -6) { 14.\x};

\foreach \x in {33,...,38}
\node at (-64+2*\x, -8) { 14.\x};

\end{scope}

\foreach \a/\b/\x/\y in {
11/13/0/0, 5/12/0/0, 0/2/0/0, 1/8/0/0, 3/10/0/0, 4/6/0/0, 7/9/0/0,
1/12/1/0, 6/13/1/0, 0/7/1/0, 2/9/1/0, 3/10/1/0, 4/11/1/0, 5/8/1/0,  
4/7/2/0, 5/11/2/0, 6/12/2/0, 0/8/2/0, 1/9/2/0, 3/10/2/0, 2/13/2/0, 
2/6/3/0, 3/10/3/0, 0/4/3/0, 5/11/3/0, 7/12/3/0, 8/13/3/0, 1/9/3/0,
1/12/4/0, 4/13/4/0, 0/9/4/0, 2/6/4/0, 3/10/4/0, 5/8/4/0, 7/11/4/0,
0/11/5/0, 8/12/5/0, 2/13/5/0, 1/5/5/0, 3/7/5/0, 4/9/5/0, 6/10/5/0, 
1/11/6/0, 2/12/6/0, 5/13/6/0, 0/8/6/0, 3/7/6/0, 4/9/6/0, 6/10/6/0, 
0/4/7/0, 1/7/7/0, 2/11/7/0, 3/8/7/0, 5/10/7/0, 6/12/7/0, 9/13/7/0,
0/2/0/1, 1/5/0/1, 3/10/0/1, 4/6/0/1, 7/9/0/1, 8/12/0/1, 11/13/0/1, 
2/4/1/1, 3/10/1/1, 5/8/1/1, 6/13/1/1, 0/7/1/1, 9/11/1/1, 1/12/1/1, 
1/4/2/1, 2/5/2/1, 3/10/2/1, 6/9/2/1, 7/12/2/1, 8/13/2/1, 0/11/2/1, 
2/13/3/1, 0/4/3/1, 1/5/3/1, 3/10/3/1, 6/9/3/1, 7/11/3/1, 8/12/3/1,
11/13/4/1, 1/12/4/1, 0/2/4/1, 3/10/4/1, 4/7/4/1, 5/8/4/1, 6/9/4/1,  
10/12/5/1, 2/11/5/1, 5/13/5/1, 0/8/5/1, 1/3/5/1, 4/7/5/1, 6/9/5/1, 
0/3/6/1, 1/7/6/1, 2/8/6/1, 4/9/6/1, 5/11/6/1, 6/12/6/1, 10/13/6/1, 
2/5/7/1, 3/9/7/1, 4/10/7/1, 6/13/7/1, 0/7/7/1, 8/11/7/1, 1/12/7/1,
1/5/0/2, 2/7/0/2, 3/10/0/2, 4/13/0/2, 6/11/0/2, 8/12/0/2, 0/9/0/2,
1/12/1/2, 4/13/1/2, 0/9/1/2, 2/7/1/2, 3/10/1/2, 5/8/1/2, 6/11/1/2,
5/7/2/2, 6/11/2/2, 8/13/2/2, 1/9/2/2, 3/10/2/2, 2/12/2/2, 0/4/2/2,
2/4/3/2, 3/10/3/2, 5/8/3/2, 6/12/3/2, 7/13/3/2, 1/9/3/2, 0/11/3/2, 
11/13/4/2, 1/12/4/2, 0/2/4/2, 3/6/4/2, 4/8/4/2, 5/9/4/2, 7/10/4/2,
0/11/5/2, 1/12/5/2, 2/13/5/2, 3/6/5/2, 4/8/5/2, 5/9/5/2, 7/10/5/2, 
0/12/6/2, 1/13/6/2, 2/11/6/2, 3/6/6/2, 4/8/6/2, 5/9/6/2, 7/10/6/2, 
0/3/7/2, 1/5/7/2, 2/6/7/2, 4/9/7/2, 7/11/7/2, 8/12/7/2, 10/13/7/2, 
11/13/0/3, 1/12/0/3, 0/2/0/3, 3/10/0/3, 4/6/0/3, 5/8/0/3, 7/9/0/3, 
0/11/1/3, 1/12/1/3, 2/13/1/3, 3/10/1/3, 4/7/1/3, 5/8/1/3, 6/9/1/3,
1/13/2/3, 0/4/2/3, 2/7/2/3, 3/10/2/3, 5/9/2/3, 6/11/2/3, 8/12/2/3,
2/13/3/3, 0/4/3/3, 1/9/3/3, 3/10/3/3, 5/8/3/3, 6/11/3/3, 7/12/3/3, 
0/3/4/3, 1/5/4/3, 2/11/4/3, 4/8/4/3, 6/10/4/3, 7/12/4/3, 9/13/4/3, 
0/3/5/3, 1/6/5/3, 2/10/5/3, 4/8/5/3, 5/11/5/3, 7/12/5/3, 9/13/5/3,  
0/3/6/3, 1/10/6/3, 2/6/6/3, 4/8/6/3, 5/11/6/3, 7/12/6/3, 9/13/6/3, 
0/3/7/3, 1/5/7/3, 2/8/7/3, 4/9/7/3, 6/11/7/3, 7/12/7/3, 10/13/7/3,
1/12/1/4, 9/13/1/4, 0/4/1/4, 2/6/1/4, 3/10/1/4, 5/8/1/4, 7/11/1/4,
5/8/2/4, 6/11/2/4, 7/12/2/4, 0/9/2/4, 3/10/2/4, 2/13/2/4, 1/4/2/4, 
1/13/3/4, 0/5/3/4, 2/8/3/4, 3/10/3/4, 4/6/3/4, 7/11/3/4, 9/12/3/4, 
2/4/4/4, 3/9/4/4, 5/11/4/4, 6/13/4/4, 0/7/4/4, 8/10/4/4, 1/12/4/4, 
0/2/5/4, 1/5/5/4, 3/9/5/4, 4/6/5/4, 7/11/5/4, 8/12/5/4, 10/13/5/4,
0/2/6/4, 1/4/6/4, 3/9/6/4, 5/8/6/4, 6/11/6/4, 7/12/6/4, 10/13/6/4
}
\draw[shift={(2*\x cm,-2*\y cm)}]
	(25.714*\a+12.857:0.87) to[out=25.714*\a+180+12.957, in=25.714*\b+180+12.857] (25.714*\b+12.857:0.87);

\foreach \a/\x/\y in {
	0/0/0, 3/0/0, 11/0/0,
	3/1/0, 5/1/0, 9/1/0, 11/1/0,
	3/2/0, 6/2/0, 11/2/0, 13/2/0,
	3/3/0, 6/3/0, 11/3/0,
	0/4/0, 4/4/0, 10/4/0,
	0/5/0, 2/5/0, 5/5/0, 9/5/0, 12/5/0,
	0/6/0, 5/6/0, 9/6/0,
	0/7/0, 3/7/0, 5/7/0, 9/7/0, 11/7/0,
	0/0/1, 3/0/1, 11/0/1,
	3/1/1, 4/1/1, 10/1/1, 11/1/1,
	1/2/1, 3/2/1, 5/2/1, 11/2/1,
	3/3/1, 8/3/1, 11/3/1, 13/3/1,
	0/4/1, 3/4/1, 11/4/1,
	0/5/1, 5/5/1, 7/5/1, 9/5/1,
	0/6/1, 4/6/1, 10/6/1,
	0/7/1, 2/7/1, 6/7/1, 8/7/1, 12/7/1,
	3/0/2, 7/0/2, 11/0/2,
	3/1/2, 7/1/2, 11/1/2,
	3/2/2, 6/2/2, 7/2/2, 11/2/2, 12/2/2,
	1/3/2, 3/3/2, 4/3/2, 10/3/2, 11/3/2,
	0/4/2, 3/4/2, 7/4/2, 11/4/2,
	0/5/2, 2/5/2, 12/5/2,
	0/6/2, 1/6/2, 2/6/2, 12/6/2, 13/6/2,
	0/7/2, 4/7/2, 10/7/2,
	0/0/3, 3/0/3, 11/0/3,
	1/1/3, 3/1/3, 11/1/3, 13/1/3,
	3/2/3, 7/2/3, 11/2/3,
	3/3/3, 7/3/3, 11/3/3, 13/3/3,
	0/4/3, 4/4/3, 9/4/3,
	0/5/3, 4/5/3, 9/5/3,
	0/6/3, 4/6/3, 9/6/3,
	0/7/3, 4/7/3, 10/7/3,
	3/1/4, 6/1/4, 8/1/4, 11/1/4,
	3/2/4, 7/2/4, 11/2/4, 13/2/4,
	3/3/4, 8/3/4, 11/3/4,
	0/4/4, 6/4/4, 8/4/4, 11/4/4,
	0/5/4, 3/5/4, 10/5/4,
	0/6/4, 3/6/4, 10/6/4
	}
\fill[shift={(2*\x cm,-2*\y cm)}]
	(25.714*\a:0.9) circle (0.05);

\foreach \a/\x/\y in {
	4/0/0, 7/0/0, 10/0/0,
	2/1/0, 4/1/0, 10/1/0, 12/1/0,
	1/2/0, 4/2/0, 8/2/0, 10/2/0,
	1/3/0, 4/3/0, 10/3/0,
	1/4/0, 5/4/0, 9/4/0, 13/4/0,
	4/5/0, 7/5/0, 10/5/0,
	4/6/0, 7/6/0, 10/6/0,
	2/7/0, 7/7/0, 12/7/0,
	4/0/1, 7/0/1, 10/0/1,
	2/1/1, 5/1/1, 9/1/1, 12/1/1,
	2/2/1, 4/2/1, 6/2/1, 10/2/1,
	1/3/1, 4/3/1, 6/3/1, 10/3/1,
	4/4/1, 6/4/1, 8/4/1, 10/4/1,
	2/5/1, 3/5/1, 11/5/1, 12/5/1,
	2/6/1, 5/6/1, 7/6/1, 9/6/1, 12/6/1,
	1/7/1, 7/7/1, 13/7/1,
	0/0/2, 4/0/2, 10/0/2,
	0/1/2, 4/1/2, 10/1/2,
	1/2/2, 4/2/2, 10/2/2,
	2/3/2, 5/3/2, 9/3/2,
	1/4/2, 2/4/2, 12/4/2, 13/4/2,
	1/5/2, 3/5/2, 7/5/2, 11/5/2, 13/5/2,
	3/6/2, 7/6/2, 11/6/2,
	2/7/2, 5/7/2, 7/7/2, 9/7/2, 12/7/2,
	4/0/3, 7/0/3, 10/0/3,
	4/1/3, 6/1/3, 8/1/3, 10/1/3,
	0/2/3, 1/2/3, 4/2/3, 5/2/3, 10/2/3,
	1/3/3, 4/3/3, 10/3/3,
	1/4/3, 3/4/3, 6/4/3, 11/4/3,
	1/5/3, 3/5/3, 10/5/3, 13/5/3,
	1/6/3, 3/6/3, 6/6/3, 11/6/3,
	1/7/3, 3/7/3, 6/7/3, 8/7/3, 12/7/3,
	1/1/4, 4/1/4, 10/1/4, 13/1/4,
	0/2/4, 2/2/4, 4/2/4, 10/2/4,
	4/3/4, 7/3/4, 10/3/4, 12/3/4,
	1/4/4, 7/4/4, 13/4/4,
	1/5/4, 2/5/4, 5/5/4, 6/5/4,
	1/6/4, 2/6/4, 4/6/4, 5/6/4, 9/6/4
	}
\fill[red, shift={(2*\x cm,-2*\y cm)}]
	(25.714*\a:0.9) circle (0.05);

\foreach \a/\x/\y in {
	1/0/0, 2/0/0, 8/0/0, 9/0/0,
	0/1/0, 6/1/0, 8/1/0,
	0/2/0, 2/2/0, 9/2/0,
	2/3/0, 7/3/0, 9/3/0, 13/3/0,
	2/4/0, 7/4/0, 12/4/0,
	1/5/0, 6/5/0, 11/5/0,
	1/6/0, 3/6/0, 8/6/0, 12/6/0,
	1/7/0, 4/7/0, 8/7/0,
	1/0/1, 2/0/1, 5/0/1, 6/0/1,
	0/1/1, 6/1/1, 8/1/1,
	7/2/1, 9/2/1, 13/2/1,
	0/3/1, 2/3/1, 5/3/1,
	5/4/1, 7/4/1, 9/4/1,
	1/5/1, 4/5/1, 8/5/1,
	1/6/1, 3/6/1, 8/6/1,
	3/7/1, 5/7/1, 10/7/1,
	2/0/2, 5/0/2, 8/0/2, 13/0/2,
	1/1/2, 5/1/2, 9/1/2, 13/1/2,
	0/2/2, 5/2/2, 8/2/2,
	6/3/2, 8/3/2, 13/3/2,
	4/4/2, 6/4/2, 9/4/2,
	4/5/2, 6/5/2, 9/5/2,
	4/6/2, 6/6/2, 9/6/2,
	1/7/2, 3/7/2, 6/7/2,
	5/0/3, 6/0/3, 8/0/3, 9/0/3,
	5/1/3, 7/1/3, 9/1/3,
	6/2/3, 9/2/3, 12/2/3,
	0/3/3, 2/3/3, 5/3/3, 9/3/3,
	2/4/3, 5/4/3, 8/4/3, 12/4/3,
	2/5/3, 6/5/3, 11/5/3,
	2/6/3, 7/6/3, 10/6/3, 13/6/3,
	2/7/3, 5/7/3, 9/7/3,
	0/1/4, 5/1/4, 9/1/4,
	1/2/4, 5/2/4, 9/2/4,
	0/3/4, 1/3/4, 5/3/4, 6/3/4,
	2/4/4, 5/4/4, 12/4/4,
	4/5/4, 7/5/4, 9/5/4, 12/5/4,
	6/6/4, 8/6/4, 12/6/4
	}
\fill[blue, shift={(2*\x cm,-2*\y cm)}]
	(25.714*\a:0.9) circle (0.05);

\foreach \a/\x/\y in {
	5/0/0, 6/0/0, 12/0/0, 13/0/0,
	1/1/0, 7/1/0, 13/1/0,
	5/2/0, 7/2/0, 12/2/0,
	0/3/0, 5/3/0, 8/3/0, 12/3/0,
	3/4/0, 6/4/0, 8/4/0, 11/4/0,
	3/5/0, 8/5/0, 13/5/0,
	2/6/0, 6/6/0, 11/6/0, 13/6/0,
	6/7/0, 10/7/0, 13/7/0,
	8/0/1, 9/0/1, 12/0/1, 13/0/1,
	1/1/1, 7/1/1, 13/1/1,
	0/2/1, 8/2/1, 12/2/1,
	7/3/1, 9/3/1, 12/3/1,
	1/4/1, 2/4/1, 12/4/1, 13/4/1,
	6/5/1, 10/5/1, 13/5/1,
	6/6/1, 11/6/1, 13/6/1,
	4/7/1, 9/7/1, 11/7/1,
	1/0/2, 6/0/2, 9/0/2, 12/0/2,
	2/1/2, 6/1/2, 8/1/2, 12/1/2,
	2/2/2, 9/2/2, 13/2/2,
	0/3/2, 7/3/2, 12/3/2,
	5/4/2, 8/4/2, 10/4/2,
	5/5/2, 8/5/2, 10/5/2,
	5/6/2, 8/6/2, 10/6/2,
	8/7/2, 11/7/2, 13/7/2,
	1/0/3, 2/0/3, 12/0/3, 13/0/3,
	0/1/3, 2/1/3, 12/1/3,
	2/2/3, 8/2/3, 13/2/3,
	6/3/3, 8/3/3, 12/3/3, 
	7/4/3, 10/4/3, 13/4/3,
	5/5/3, 8/5/3, 12/5/3,
	5/6/3, 8/6/3, 12/6/3,
	7/7/3, 11/7/3, 13/7/3,
	2/1/4, 7/1/4, 12/1/4,
	6/2/4, 8/2/4, 12/2/4,
	2/3/4, 9/3/4, 13/3/4,
	3/4/4, 4/4/4, 9/4/4, 10/4/4,
	8/5/4, 11/5/4, 13/5/4,
	7/6/4, 11/6/4, 13/6/4
	}
\fill[green, shift={(2*\x cm,-2*\y cm)}]
	(25.714*\a:0.9) circle (0.05);

\end{scope}

\end{tikzpicture}
\caption{Tilings of $2{\bb T}^2$ by single 8-gon, 10-gon, 12-gon, 14-gon.} 
\label{2T2_n8}
\end{figure}

\begin{figure}[htp]
\centering
\begin{tikzpicture}[>=latex,scale=1]

\begin{scope}[shift={(2cm,6cm)}]

\foreach \a in {0,...,15}
\foreach \x in {0,...,5}
\foreach \y in {0,1}
\draw[shift={(2*\x cm, -2*\y cm)}]
	(22.5*\a-11.25:0.9) -- (22.5*\a+11.25:0.9);

\foreach \a in {0,...,15}
\foreach \x in {-1,...,6}
\draw[shift={(2*\x cm, -2*2 cm)}]
	(22.5*\a-11.25:0.9) -- (22.5*\a+11.25:0.9);		
\foreach \x in {1,...,6}
\node[gray!70] at (-2+2*\x,0) {16.\x};

\foreach \x in {7,...,12}
\node[gray!70] at (-14+2*\x,-2) {16.\x};

\foreach \x in {13,...,20}
\node[gray!70] at (-28+2*\x,-4) {16.\x};

\foreach \a/\b/\x/\y in {
	1/4/0/0, 2/10/0/0, 3/11/0/0, 5/8/0/0, 6/14/0/0, 7/15/0/0, 9/12/0/0, 0/13/0/0, 
	3/15/1/0, 0/5/1/0, 1/12/1/0, 2/6/1/0, 4/9/1/0, 7/11/1/0, 8/13/1/0, 10/14/1/0,
	2/15/2/0, 0/4/2/0, 1/5/2/0, 3/6/2/0, 7/10/2/0, 8/12/2/0, 9/13/2/0, 11/14/2/0, 
	0/3/3/0, 1/4/3/0, 2/5/3/0, 6/15/3/0, 7/10/3/0, 8/12/3/0, 9/13/3/0, 11/14/3/0,
	0/2/4/0, 1/4/4/0, 3/5/4/0, 6/15/4/0, 7/10/4/0, 8/12/4/0, 9/13/4/0, 11/14/4/0, 
	4/8/5/0, 5/11/5/0, 2/6/5/0, 7/12/5/0, 9/14/5/0, 0/10/5/0, 1/13/5/0, 3/15/5/0, 
	0/2/0/1, 1/5/0/1, 3/11/0/1, 4/6/0/1, 7/10/0/1, 8/13/0/1, 9/14/0/1, 12/15/0/1, 
	1/15/1/1, 0/6/1/1, 2/9/1/1, 3/11/1/1, 4/13/1/1, 5/7/1/1, 8/12/1/1, 10/14/1/1, 
	4/6/2/1, 5/13/2/1, 7/10/2/1, 0/8/2/1, 1/9/2/1, 3/11/2/1, 12/14/2/1, 2/15/2/1, 
	1/14/3/1, 5/15/3/1, 0/6/3/1, 2/9/3/1, 3/12/3/1, 4/7/3/1, 8/11/3/1, 10/13/3/1, 
	1/15/4/1, 0/5/4/1, 2/9/4/1, 3/12/4/1, 4/6/4/1, 7/11/4/1, 8/13/4/1, 10/14/4/1, 
	3/7/5/1, 4/11/5/1, 0/5/5/1, 6/12/5/1, 8/13/5/1, 9/15/5/1, 1/10/5/1, 2/14/5/1, 
	1/14/-1/2, 4/15/-1/2, 0/10/-1/2, 2/7/-1/2, 3/11/-1/2, 5/9/-1/2, 6/12/-1/2, 8/13/-1/2,    
	1/14/0/2, 10/15/0/2, 0/4/0/2, 2/7/0/2, 3/11/0/2, 5/9/0/2, 6/12/0/2, 8/13/0/2,
	1/14/1/2, 4/15/1/2, 0/5/1/2, 2/8/1/2, 3/11/1/2, 6/10/1/2, 7/12/1/2, 9/13/1/2,
	1/4/2/2, 2/6/2/2, 3/11/2/2, 5/12/2/2, 7/10/2/2, 8/14/2/2, 9/15/2/2, 0/13/2/2,
	1/14/3/2, 6/15/3/2, 0/7/3/2, 2/9/3/2, 3/11/3/2, 4/12/3/2, 5/8/3/2, 10/13/3/2, 
	3/11/4/2, 4/13/4/2, 5/15/4/2, 6/9/4/2, 1/7/4/2, 2/8/4/2, 10/14/4/2, 0/12/4/2, 
	2/15/5/2, 0/4/5/2, 1/5/5/2, 3/10/5/2, 6/9/5/2, 7/12/5/2, 8/13/5/2, 11/14/5/2, 
	1/5/6/2, 2/7/6/2, 3/10/6/2, 4/14/6/2, 6/11/6/2, 8/13/6/2, 9/15/6/2, 0/12/6/2
	}
\draw[shift={(2*\x cm, -2*\y cm)}]
	(22.5*\a+22.5:0.88) to[out=22.5*\a+180+22.5, in=22.5*\b+180+22.5] (22.5*\b+22.5:0.88);

\foreach \a/\x/\y in { 
	3/0/0, 10/0/0, 12/0/0,
	2/1/0, 7/1/0, 12/1/0,
	0/2/0, 2/2/0, 5/2/0,
	1/3/0, 3/3/0, 5/3/0,
	0/4/0, 3/4/0, 6/4/0,
	2/5/0, 7/5/0, 13/5/0,
	0/0/1, 3/0/1, 12/0/1,
	3/1/1, 9/1/1, 12/1/1,
	3/2/1, 12/2/1, 15/2/1,
	2/3/1, 10/3/1, 14/3/1,
	3/4/1, 9/4/1, 13/4/1,
	3/5/1, 8/5/1, 14/5/1,
	3/-1/2, 7/-1/2, 12/-1/2,
	3/0/2, 7/0/2, 12/0/2,
	3/1/2, 8/1/2, 12/1/2,
	3/2/2, 6/2/2, 12/2/2,
	3/3/2, 5/3/2, 9/3/2, 12/3/2,
	1/4/2, 3/4/2, 8/4/2, 12/4/2,
	15/5/2, 3/5/2, 11/5/2,
	3/6/2, 7/6/2, 11/6/2
	}
\fill[shift={(2*\x cm, -2*\y cm)}]
	(22.5*\a+11.75:0.9) circle (0.05);

\foreach \a/\x/\y in { 
	2/0/0, 4/0/0, 11/0/0,
	4/1/0, 10/1/0, 15/1/0,
	1/2/0, 4/2/0, 6/2/0,
	7/3/0, 11/3/0, 15/3/0,
	7/4/0, 11/4/0, 15/4/0,
	4/5/0, 9/5/0, 15/5/0,
	4/0/1, 7/0/1, 11/0/1,
	4/1/1, 11/1/1, 14/1/1,
	4/2/1, 7/2/1, 11/2/1,
	4/3/1, 8/3/1, 12/3/1,
	8/4/1, 11/4/1, 14/4/1,
	1/5/1, 5/5/1, 11/5/1,
	0/-1/2, 4/-1/2, 11/-1/2,
	1/0/2, 4/0/2, 11/0/2, 15/0/2,
	0/1/2, 4/1/2, 6/1/2, 11/1/2,
	2/2/2, 4/2/2, 7/2/2, 11/2/2,
	4/3/2, 11/3/2, 13/3/2,
	4/4/2, 11/4/2, 14/4/2,
	1/5/2, 4/5/2, 6/5/2, 10/5/2,
	4/6/2, 10/6/2, 15/6/2
	}
\fill[gray,shift={(2*\x cm, -2*\y cm)}]
	(22.5*\a+11.75:0.9) circle (0.05);

\foreach \a/\x/\y in { 
	1/0/0, 5/0/0, 9/0/0, 13/0/0,
	1/1/0, 5/1/0, 9/1/0, 13/1/0,
	3/2/0, 7/2/0, 11/2/0, 15/2/0,
	0/3/0, 2/3/0, 4/3/0, 6/3/0,
	1/4/0, 2/4/0, 4/4/0, 5/4/0,
	0/5/0, 3/5/0, 6/5/0, 11/5/0,
	1/0/1, 2/0/1, 5/0/1, 6/0/1,
	0/1/1, 1/1/1, 6/1/1, 7/1/1,
	5/2/1, 6/2/1, 13/2/1, 14/2/1,
	3/3/1, 9/3/1, 11/3/1, 13/3/1,
	0/4/1, 1/4/1, 5/4/1, 6/4/1,
	0/5/1, 6/5/1, 9/5/1, 13/5/1,
	1/-1/2, 5/-1/2, 10/-1/2, 15/-1/2,
	2/0/2, 8/0/2, 14/0/2,
	2/1/2, 9/1/2, 14/1/2,
	1/2/2, 5/2/2, 13/2/2,
	2/3/2, 10/3/2, 14/3/2,
	2/4/2, 7/4/2, 9/4/2,
	0/5/2, 2/5/2, 5/5/2,
	2/6/2, 5/6/2, 8/6/2, 14/6/2
	 }
\fill[red,shift={(2*\x cm, -2*\y cm)}]
	(22.5*\a+11.75:0.9) circle (0.05);

\foreach \a/\x/\y in { 
	0/0/0, 7/0/0, 14/0/0,
	0/1/0, 3/1/0, 6/1/0,
	9/2/0, 12/2/0, 14/2/0,
	8/3/0, 10/3/0, 13/3/0,
	8/4/0, 10/4/0, 13/4/0,
	1/5/0, 10/5/0, 14/5/0,
	8/0/1, 10/0/1, 14/0/1,
	2/1/1, 10/1/1, 15/1/1,
	1/2/1, 8/2/1, 10/2/1,
	0/3/1, 5/3/1, 7/3/1,
	2/4/1, 10/4/1, 15/4/1,
	2/5/1, 10/5/1, 15/5/1,
	2/-1/2, 8/-1/2, 14/-1/2,
	0/0/2, 5/0/2, 10/0/2,
	1/1/2, 5/1/2, 15/1/2,
	8/2/2, 10/2/2, 15/2/2,
	0/3/2, 6/3/2, 8/3/2,
	0/4/2, 5/4/2, 13/4/2,
	13/5/2, 7/5/2, 9/5/2,
	1/6/2, 6/6/2, 12/6/2
	 }
\fill[blue,shift={(2*\x cm, -2*\y cm)}]
	(22.5*\a+11.75:0.9) circle (0.05);

\foreach \a/\x/\y in { 
	6/0/0, 8/0/0, 15/0/0,
	8/1/0, 11/1/0, 14/1/0,
	8/2/0, 10/2/0, 13/2/0,
	9/3/0, 12/3/0, 14/3/0,
	9/4/0, 12/4/0, 14/4/0,
	5/5/0, 8/5/0, 12/5/0,
	9/0/1, 13/0/1, 15/0/1,
	5/1/1, 8/1/1, 13/1/1,
	0/2/1, 2/2/1, 9/2/1,
	1/3/1, 6/3/1, 15/3/1,
	4/4/1, 7/4/1, 12/4/1,
	4/5/1, 7/5/1, 12/5/1,
	6/-1/2, 9/-1/2, 13/-1/2,
	6/0/2, 9/0/2, 13/0/2,
	7/1/2, 10/1/2, 13/1/2,
	0/2/2, 9/2/2, 14/2/2,
	1/3/2, 7/3/2, 15/3/2,
	6/4/2, 10/4/2, 15/4/2,
	8/5/2, 12/5/2, 14/5/2,
	0/6/2, 9/6/2, 13/6/2
	 }
\fill[green,shift={(2*\x cm, -2*\y cm)}]
	(22.5*\a+11.75:0.9) circle (0.05);
		
\end{scope}

\foreach \x in {0,...,7}
\foreach \a in {0,...,17}
\draw[xshift=2*\x cm]
	(20*\a:0.9) -- (20*\a+20:0.9);

\foreach \x in {1,...,8}
\node[gray!70] at (-2+2*\x,0) {18.\x};
	
\foreach \a/\b/\x in {
	0/3/0, 5/8/0, 9/12/0, 14/17/0, 1/10/0, 2/11/0, 4/13/0, 6/15/0, 7/16/0,
	0/3/1, 5/8/1, 9/12/1, 14/17/1, 1/6/1, 2/7/1, 10/15/1, 11/16/1, 4/13/1, 
	3/16/6, 6/17/6, 0/11/6, 1/14/6, 2/7/6, 4/9/6, 5/12/6, 8/13/6, 10/15/6, 
	4/7/7, 10/13/7, 0/14/7, 3/17/7, 2/15/7, 5/11/7, 6/12/7, 8/16/7, 1/9/7, 
	0/3/4, 5/8/4, 9/12/4, 14/17/4, 1/15/4, 2/16/4, 6/10/4, 7/11/4, 4/13/4,  
	2/17/3, 6/9/3, 10/14/3, 12/16/3, 0/7/3, 1/8/3, 3/11/3, 5/15/3, 4/13/3, 
	2/16/2, 11/17/2, 0/5/2, 1/12/2, 3/8/2, 4/13/2, 6/10/2, 7/14/2, 9/15/2,
	2/16/5, 7/11/5, 1/6/5, 10/15/5, 5/17/5, 0/12/5, 3/9/5, 8/14/5, 4/13/5
	}
\draw[xshift=2*\x cm]
	(20*\a+10:0.88) to[out=20*\a+190, in=20*\b+190] (20*\b+10:0.88);

\foreach \a/\x in {
	0/0, 4/0, 14/0,
	0/1, 4/1, 14/1,
	0/6, 6/6, 12/6,
	0/7, 3/7, 15/7,
	0/4, 4/4, 14/4,
	4/3, 11/3, 14/3, 
	4/2, 8/2, 14/2,
	4/5, 9/5, 14/5
	}
\fill[xshift=2*\x cm]
	(20*\a:0.9) circle (0.05);

\foreach \a/\x in {
	5/0, 13/0, 9/0,
	5/1, 9/1, 13/1,
	2/6, 8/6, 14/6,
	2/7, 9/7, 16/7,
	5/4, 9/4, 13/4,
	5/3, 13/3, 16/3,
	1/2, 5/2, 13/2,
	0/5, 5/5, 13/5
	}
\fill[gray, xshift=2*\x cm]
	(20*\a:0.9) circle (0.05);

\foreach \a/\x in {
	1/0, 11/0, 3/0,
	1/1, 3/1, 7/1,
	4/6, 10/6, 16/6,
	1/7, 10/7, 14/7,
	1/4, 16/4, 3/4,
	3/3, 12/3, 17/3,
	3/2, 9/2, 16/2,
	3/5, 10/5, 16/5
	}
\fill[red, xshift=2*\x cm]
	(20*\a:0.9) circle (0.05);

\foreach \a/\x in {
	6/0, 8/0, 16/0,
	10/1, 12/1, 16/1,
	3/6, 7/6, 17/6,
	4/7, 8/7, 17/7,
	6/4, 8/4, 11/4,
	6/3, 10/3, 15/3,
	0/2, 6/2, 11/2,
	1/5, 7/5, 12/5
	}
\fill[blue, xshift=2*\x cm]
	(20*\a:0.9) circle (0.05);

\foreach \a/\x in {
	2/0, 12/0, 10/0,
	2/1, 6/1, 8/1,
	9/6, 5/6, 13/6,
	5/7, 7/7, 12/7,
	2/4, 15/4, 17/4,
	0/3, 2/3, 8/3,
	2/2, 12/2, 17/2,
	2/5, 6/5, 17/5
	}
\fill[green, xshift=2*\x cm]
	(20*\a:0.9) circle (0.05);

\foreach \a/\x in {
	7/0, 15/0, 17/0,
	11/1, 15/1, 17/1,
	1/6, 11/6, 15/6,
	6/7, 11/7, 13/7,
	7/4, 10/4, 12/4,
	1/3, 7/3, 9/3,
	7/2, 10/2, 15/2,
	8/5, 11/5, 15/5
	}
\fill[brown, xshift=2*\x cm]
	(20*\a:0.9) circle (0.05);

\end{tikzpicture}
\caption{Tilings of $2{\bb T}^2$ by single 16-gon, 18-gon.} 
\label{2T2_n16}
\end{figure}

Zamorzaeva-Orleanschi \cite{zamo1,zamo2} also studied single tile tilings of $2{\bb T}^2$ and came up with a list. She missed the following tilings: 12.15, 12.19, 12.21, 12.22, 12.26, 14.34, 14.36, 18.3, 18.6. Her list also duplicated the following tilings: 12:20, 12.24, 14.3, 14.12, 14.28, 16.18.

\section{Single Tile Tiling of Non-orientable Surface}

Both opposing and twisted edge pairs in Figure \ref{edge_pair2} may appear in the planar diagram of a single tile tiling of a non-orientable surface. In the planar diagrams and the corresponding vertex sets, we need to keep track of the orientation. 

We denote an opposing edge pair by $\xar{i}\xar{j}_{\sigma}$, where $\sigma=+$ means the pair is opposing, and $\sigma=-$ means the pair is twisted. There is no order in edge pair. In other words, the edge pairs $\xar{i}\xar{j}_{\sigma}$ and $\xar{j}\xar{i}_{\sigma}$ are the same. The non-orientability of the surface means there is at least one twisted pair in the planar diagram. Figure \ref{edge_pair3} is a decagon planar diagram $D=(\xar{0}\xar{1}_-,\xar{2}\xar{5}_+,\xar{3}\xar{8}_-,\xar{4}\xar{7}_-,\xar{6}\xar{9}_+)$. The solid chords mean opposing edge pairs, and the dashed chords mean twisted edge pairs. The usual notation for the planar diagram is $a_1a_1a_2a_3a_4a_2^{-1}a_5a_4a_3a_5^{-1}$.

\begin{figure}[htp]
\centering
\begin{tikzpicture}[>=latex,scale=1]

\foreach \a in {0,...,9}
{
\foreach \b in {0,1}
\draw[xshift=4*\b cm]
	(36*\a:1.2) -- (36*\a+36:1.2);
	
\node at (36*\a:1) {\footnotesize \a};

}

\foreach \a/\b in {
	2/5, 6/9}
\draw
	(36*\a+18:1.14) to[out=36*\a+198, in=36*\b+198] (36*\b+18:1.14);

\foreach \a/\b in {
	0/1, 3/8, 4/7}
\draw[dash pattern=on 2pt off 1pt]
	(36*\a+18:1.14) to[out=36*\a+198, in=36*\b+198] (36*\b+18:1.14);
	
\node at (0,-1.8) {\footnotesize $(\xar{0}\xar{1}_-,\xar{2}\xar{5}_+,\xar{3}\xar{8}_-,\xar{4}\xar{7}_-,\xar{6}\xar{9}_+)$};

\begin{scope}[xshift=4cm]

\foreach \a/\b in {
	0/1, 1/1, 2/1, 5/-1, 3/1, 8/1, 4/1, 7/1, 6/1, 9/-1}
\draw[->]
	(36*\a+18:1.14) -- ++(36*\a+108:0.1*\b);

\foreach \a/\b in {
	0/1, 1/1, 2/2, 5/2, 3/3, 8/3, 4/4, 7/4, 6/5, 9/5}
\node at (36*\a+18:1.35) {\footnotesize $a_{\b}$};	
	
\node at (0,-1.8) {\footnotesize $a_1a_1a_2a_3a_4a_2^{-1}a_5a_4a_3a_5^{-1}$};

\end{scope}

\end{tikzpicture}
\caption{Planar diagram for a non-orientable surface.} 
\label{edge_pair3}
\end{figure}

Figure \ref{vertex_twist} shows how a twisted pair of edges are glued together. We decorate all the corners and edges by the orientation $\sigma=\pm$ in the picture. We split the twisted edge pair $\xar{i}\xar{j}_-$ into ordered pairs $\xar{i}_-\xar{j}_+$ and $\xar{j}_+\xar{i}_-$, and get
\[
\xar{i}_-\xar{j}_+
\implies \bullet=i_-j_+\cdots, \quad
\xar{j}_+\xar{i}_-
\implies \circ=(j+1)_+(i+1)_-\cdots.
\]

\begin{figure}[h]
\centering
\begin{tikzpicture}[>=latex,scale=1]

\draw
	(0,-1.2) -- (0,1.2)
	(-2,-1.2) -- (2,-1.2)
	(-2,1.2) -- (2,1.2);

\fill 
	(0,-1.2) circle (0.08);
\filldraw[fill=white]
	(0,1.2) circle (0.08);
		
\draw[yshift=1.2 cm, <-]
	(15:1.05) arc (15:-195:1.05);
	
\draw[yshift=-1.2 cm, ->]
	(-15:0.6) arc (-15:195:0.6);

\draw[xshift=-1.5cm, ->]
	(120:0.3) arc (120:420:0.3);
\draw[xshift=1.5cm, <-]
	(120:0.3) arc (120:420:0.3);

\node at (-1.5,0) {$+$};
\node at (1.5,0) {$-$};

\node at (0.25,-0.2) {\footnotesize $\xar{i}_-$};
\node at (-0.25,-0.2) {\footnotesize $\xar{j}_+$};

\node at (1.7,0.9) {\footnotesize $\overline{(i\!+\!1)}_-$};
\node at (-1.6,0.9) {\footnotesize $\overline{(j\!+\!1)}_+$};
\node at (1.6,-0.9) {\footnotesize $\overline{(i\!-\!1)}_-$};
\node at (-1.6,-0.9) {\footnotesize $\overline{(j\!-\!1)}_+$};

\node at (0.2,-1) {\scriptsize $i_-$};
\node at (-0.25,-1) {\scriptsize $j_+$};
\node at (0.5,1) {\scriptsize $(i\!+\!1)_-$};
\node at (-0.5,1) {\scriptsize $(j\!+\!1)_+$};

\end{tikzpicture}
\caption{Twisted edge pair induces adjacent corners at vertices.} 
\label{vertex_twist}
\end{figure}

The glueing of an opposing pair is still described by Figure \ref{vertex_oppose}, in which all corners and edges are decorated by $+$. Then the picture means that we split an opposing pair $\xar{i}\xar{j}_+$ into ordered pairs $\xar{i}_+\xar{j}_+$ and $\xar{j}_+\xar{i}_+$, and get
\[
\xar{i}_+\xar{j}_+
\implies \bullet =(i+1)_+j_+\cdots, \quad
\xar{j}_+\xar{i}_+
\implies \circ=(j+1)_+i_+\cdots.
\]
The glueing of an opposing pair is equivalently described by the flip of Figure \ref{vertex_oppose}, with all corners and edges decorated by $-$. This means that we split $\xar{i}\xar{j}_+$ into ordered pairs $\xar{i}_-\xar{j}_-$ and $\xar{j}_-\xar{i}_-$, and get
\[
\xar{i}_-\xar{j}_-
\implies \bullet=i_-(j+1)_-\cdots,\quad
\xar{j}_-\xar{i}_-
\implies \circ=j_-(i+1)_-\cdots.
\]

The adjacent corners combine to give vertices. A vertex is a circularly ordered subset of ${\bb Z}_n$, with each corner decorated by $\sigma=\pm$. By circularly ordered, we mean rotation (of circular order) and preserving all $\sigma$, or reversion (of circular order) and changing all $\sigma$ to $-\sigma$, give the same vertex. For example, we have
\[
0_+6_+2_+1_-
\overset{\text{rotation}}{=\joinrel=\joinrel=\joinrel=}
2_+1_-0_+6_+
\overset{\text{reversion}}{=\joinrel=\joinrel=\joinrel=}
6_-0_-1_+2_-
\overset{\text{rotation}}{=\joinrel=\joinrel=\joinrel=}
0_-1_+2_-6_-.
\]

We may summarise the formulae from edge pairs to adjacent corners as the following
\begin{align*}
i_+\cdots
&\overset{\scriptsize \xar{i-1}_+\xar{j}_{\sigma}}{=\joinrel=\joinrel=\joinrel=}
\begin{cases}
i_+j_+\cdots, &\text{if }\sigma=+ \\
i_+(j+1)_-\cdots, &\text{if }\sigma=-
\end{cases}, \\
i_-\cdots
&\overset{\scriptsize \xar{i}_-\xar{j}_{\sigma}}{=\joinrel=\joinrel=\joinrel=}
\begin{cases}
i_-j_+\cdots, &\text{if }\sigma=+ \\
i_-(j+1)_-\cdots, &\text{if }\sigma=-
\end{cases}.
\end{align*}
For example, the planar diagram $D=(\xar{0}\xar{1}_-,\xar{2}\xar{5}_+,\xar{3}\xar{8}_-,\xar{4}\xar{7}_-,\xar{6}\xar{9}_+)$ in Figure \ref{edge_pair3} implies the vertices 
\begin{align*}
0_+\cdots
&\overset{\scriptsize \xar{9}_+\xar{6}_+}{=\joinrel=\joinrel=}
0_+6_+\cdots
\overset{\scriptsize \xar{5}_+\xar{2}_+}{=\joinrel=\joinrel=}
0_+6_+2_+\cdots
\overset{\scriptsize \xar{1}_+\xar{0}_-}{=\joinrel=\joinrel=}
0_+6_+2_+1_-\cdots
\overset{\scriptsize \xar{1}_-\xar{0}_+}{=\joinrel=\joinrel=}
0_+6_+2_+1_-, \\
3_+\cdots 
&\overset{\scriptsize \xar{2}_+\xar{5}_+}{=\joinrel=\joinrel=} 
3_+5_+\cdots
\overset{\scriptsize \xar{4}_+\xar{7}_-}{=\joinrel=\joinrel=} 
3_+5_+8_-\cdots
\overset{\scriptsize \xar{8}_-\xar{3}_+}{=\joinrel=\joinrel=} 
3_+5_+8_-, \\
4_+\cdots
&\overset{\scriptsize \xar{3}_+\xar{8}_-}{=\joinrel=\joinrel=}  
4_+9_-\cdots
\overset{\scriptsize \xar{9}_-\xar{6}_-}{=\joinrel=\joinrel=}  
4_+9_-7_-\cdots
\overset{\scriptsize \xar{7}_-\xar{4}_+}{=\joinrel=\joinrel=}  
4_+9_-7_-.
\end{align*}
We get the vertex set $(0_+6_+2_+1_-,3_+5_+8_-,4_+9_-7_-)$ at the end.

We can easily read vertices directly from the picture of the planar diagram, with extra attention to the direction in determining the signs of the corners. If in reading the vertex, we arrive at a corner $i_+$ with positive sign, then we look at the edge before $i$ (i.e., on the clockwise side of the corner). If the edge is in an opposing pair (the chord is solid), then we go to the next corner by following the rule on the right of Figure \ref{pd_vs}, and the next corner still has positive sign. If the edge is in a twisted pair (the chord is dashed), then we go to the next corner by following the rule on the right of Figure \ref{pd_vs2}, and the sign of the next corner is changed to negative. 

\begin{figure}[htp]
\centering
\begin{tikzpicture}[>=latex,scale=1]

\draw[gray, very thick]
	(36*3:0.9) to[out=180+36*3, in=180+36*5] (36*5:0.9) to[out=180+36*5, in=180+36*8] (36*8:0.9) to[out=180+36*8, in=180+36*3] (36*3:0.9);

\draw[gray, very thick, ->]
	(0,0) -- ++(108:0.3);
	
\foreach \a in {0,...,9}
\draw
	(36*\a:1.2) -- (36*\a+36:1.2);
	
\foreach \a/\b in {
	2/5, 6/9}
\draw
	(36*\a+18:1.14) to[out=36*\a+198, in=36*\b+198] (36*\b+18:1.14);

\foreach \a/\b in {
	0/1, 3/8, 4/7}
\draw[dash pattern=on 2pt off 1pt]
	(36*\a+18:1.14) to[out=36*\a+198, in=36*\b+198] (36*\b+18:1.14);

\node at (108:1.05) {\scriptsize $+$};
\node at (180:1.05) {\scriptsize $+$};
\node at (-72:1.05) {\scriptsize $-$};

\begin{scope}[xshift=4cm]

\draw
	(50:1.2) -- (80:1.2)
	(150:1.2) -- (180:1.2);

\draw[dashed]
	(80:1.2) arc (80:150:1.2)
	(180:1.2) arc (180:410:1.2);

\draw[dash pattern=on 2pt off 1pt]
	(65:1.16) to[out=245, in=-15] (165:1.16);

\draw[gray, very thick, ->]
	(50:0.9) to[out=180+50, in=180+150] (150:0.95);

\draw[gray, very thick, ->]
	(180:0.9) to[out=0, in=180+80] (80:0.95);

	
\fill
	(50:1.2) circle (0.05)
	(150:1.2) circle (0.05);
		
\filldraw[fill=white]
	(80:1.2) circle (0.05)
	(180:1.2) circle (0.05);

\node at (180:1) {\scriptsize $+$};
\node at (50:1) {\scriptsize $-$};
\node at (150:1) {\scriptsize $+$};
\node at (80:1) {\scriptsize $-$};
	
\end{scope}
	
\end{tikzpicture}
\caption{Planar diagram to vertices along twisted pair.} 
\label{pd_vs2}
\end{figure}

On the other hand, if we start at a corner $j_-$ with negative sign, then we look at the edge after $j$ (i.e., on the counterclockwise side of the corner). Then depending on whether the edge is in opposing or twisted pair, we proceed in similar way. 

The surface given by the planar diagram in Figure \ref{edge_pair3} has three vertices. The Euler number of the surface is $3-5+1=-1$. Therefore the surface $S_D=3{\bb P}^2$.

Conversely, we have the formulae from vertices to edge pairs
\begin{align*}
i_+j_+\cdots &\implies (\xar{i-1})\xar{j}_+, &
i_-j_+\cdots &\implies \xar{i}\xar{j}_-, \\
i_+j_-\cdots &\implies (\xar{i-1})(\xar{j-1})_-, &
i_-j_-\cdots &\implies \xar{i}(\xar{j-1})_+.
\end{align*}
For example, the vertex $0_+6_+2_+1_-$ implies edge pairs
\begin{align*}
0_+6_+\cdots &\implies \xar{9}\xar{6}_+=\xar{6}\xar{9}_+, &
6_+2_+\cdots &\implies \xar{5}\xar{2}_+=\xar{2}\xar{5}_+, \\
2_+1_-\cdots &\implies \xar{1}\xar{0}_-=\xar{0}\xar{1}_-, &
1_-0_+\cdots &\implies \xar{1}\xar{0}_-=\xar{0}\xar{1}_-.
\end{align*}
Since the edge pair $\xar{0}\xar{1}_-$ already appears twice, it will not appear from adjacent corners at the other vertices. On the other hand, the vertex set should also contain the vertices $3_+5_+\cdots$ and $7_+9_+\cdots$ implied by $\xar{2}_+\xar{5}_+$ and $\xar{6}_+\xar{9}_+$. Indeed, $3_+5_+\cdots$ is part of $3_+5_+8_-$, and $7_+9_+\cdots$ is part of $4_+9_-7_-=4_-7_+9_+$.

Similar to the orientable case, planar diagrams and vertex pairs are equivalent. The following is the analogue of Proposition \ref{avs}. The proof is completely similar and is omitted.

\begin{proposition}\label{avs2}
A disjoint union of circularly ordered subsets of ${\bb Z}_n$, together with decorations by $\sigma=\pm$, is a vertex set if and only if the following are satisfied:
\begin{itemize}
\item Each (circularly ordered) subset has at least three corners.
\item $i_+(i-1)_+\cdots$ is not a subset.
\item If $i_+j_+\cdots$ is a subset, then $(j+1)_+(i-1)_+\cdots$ is a subset. 
\item If $i_+j_-\cdots$ is a subset, then $(j-1)_-(i-1)_+\cdots$ is a subset.
\end{itemize}
\end{proposition}

The condition $i_+(i-1)_+\cdots$ is not a vertex is the same as the condition $(i-1)_-i_-\cdots$ is not a vertex. Similar remark applies to the third condition. 

Correspondingly, the following are the conditions for edge pairs in a planar diagram
\begin{itemize}
\item No degree one vertex: $\xar{i}\xar{i+1}_+$ is not an edge pair.
\item No degree two vertex in opposing way: If $\xar{i}\xar{j}_+$ is an edge pair, then $\xar{i+1}\xar{j-1}_+$ is not an edge pair.
\item No degree two vertex in twisted way: If $\xar{i}\xar{j}_-$ is an edge pair, then $\xar{i+1}\xar{j+1}_-$ is not an edge pair.
\end{itemize} 

Two planar diagrams represent the same tiling if they are related by the rotation $i\mapsto c+i$ and $\bar{i}\mapsto \overline{c+i}$, or the reversion $i\mapsto c-i$ and $\bar{i}\mapsto \overline{c-i-1}$. We note that the decorations $\sigma$ are not changed. Moreover, after applying the relabelling $c-i$ to corners, the circular order of vertices needs to be reversed. For example, the decagon planar diagram in Figure \ref{edge_pair3} and the following represent the same tiling
\begin{align*} 
(\xar{0}\xar{1}_-,\xar{2}\xar{5}_+,\xar{3}\xar{8}_-,\xar{4}\xar{7}_-,\xar{6}\xar{9}_+)
&\xrightarrow{\overline{4+i}}
(\xar{4}\xar{5}_-,\xar{6}\xar{9}_+,\xar{7}\xar{2}_-,\xar{8}\xar{1}_-,\xar{0}\xar{3}_+) \\
&=(\xar{0}\xar{3}_+,\xar{1}\xar{8}_-,\xar{2}\xar{7}_-,\xar{4}\xar{5}_-,\xar{6}\xar{9}_+), \\
(\xar{0}\xar{1}_-,\xar{2}\xar{5}_+,\xar{3}\xar{8}_-,\xar{4}\xar{7}_-,\xar{6}\xar{9}_+)
&\xrightarrow{\overline{4-i-1}}
(\xar{3}\xar{2}_-,\xar{1}\xar{8}_+,\xar{0}\xar{5}_-,\xar{9}\xar{6}_-,\xar{7}\xar{4}_+) \\
&=(\xar{0}\xar{5}_-,\xar{1}\xar{8}_+,\xar{2}\xar{3}_-,\xar{4}\xar{7}_+,\xar{6}\xar{9}_-).
\end{align*}
These are equivalent to the following changes in terms of the vertex sets
\begin{align*}
(0_+6_+2_+1_-,3_+5_+8_-,4_+9_-7_-) 
&\xrightarrow{4+i} 
(4_+0_+6_+5_-,7_+9_+2_-,8_+3_-1_-)  \\
&=(0_+6_+5_-4_+,1_-8_+3_-,2_-7_+9_+) \\
&=(0_+6_+5_-4_+,1_+3_+8_-,2_+9_-7_-). \\
(0_+6_+2_+1_-,3_+5_+8_-,4_+9_-7_-) 
&\xrightarrow{4-i} 
(4_+8_+2_+3_-,1_+9_+6_-,0_+5_-7_-)  \\
&\xrightarrow{\text{reverse}} 
(3_-2_+8_+4_+,6_-9_+1_+,7_-5_-0_+) \\
&=(0_+7_-5_-,1_+6_-9_+,2_+8_+4_+3_-).
\end{align*}

We may create computer program to find all the single tile tilings of a fixed non-orientable surface of Euler characteristic $<0$, by congruent $n$-gons with $n\ge 7$. Table \ref{P2tiling} shows the numbers of tilings for non-orientable surfaces of small genus. Figure \ref{3P2_n8} gives all the single tile tilings of $3{\bb P}^2$ by an $n$-gon with $n\ge 7$.

\begin{table}[htp]
	\centering
\begin{tabular}{|c||c|c|c|c|c|c|c|c|c|}
\hline
polygon size & 8 & 10 & 12 & 14 & 16 & 18 & 20 & 22 & 24 \\
\hline
$3{\bb P}^2$ & 22 & 24 & 11 &&&&& & \\
\hline
$4{\bb P}^2$ & 47 & 279 & 682 & 838 & 508 & 144 &&& \\
\hline
$5{\bb P}^2$ & & 473 & 4928 & 20979 & 47462 & 62283 & 47825 & 19971 & 3627 \\
\hline
\end{tabular}
\caption{Number of single tile tilings of non-orientable hyperbolic surfaces of small genus. }
\label{P2tiling}
\end{table}

\begin{figure}[htp]
\centering
\begin{tikzpicture}[>=latex,scale=1]


\begin{scope}[gray!50]

\foreach \a in {0,1,2}
\draw
	(-0.9,0.9-1.8*\a) -- ++(10.8,0);
	
\draw
	(-0.9,0.9) -- ++(0,-5.4) -- ++(5.4,0) -- ++(0,1.8)
	(4.5,0.9) -- ++(0,-1.8)
	(9.9,0.9) -- ++(0,-3.6)
	(2.7,-0.9) -- ++(0,-1.8)
	(6.3,-0.9) -- ++(0,-1.8) ;

\end{scope}

\foreach \x in {1,...,6}
\node[gray] at (-1.8+1.8*\x, 0) { 8.\x};

\foreach \x in {7,...,12}
\node[gray] at (-1.8*7+1.8*\x, -1.8) { 8.\x};

\foreach \x in {13,...,18}
\node[gray] at (-1.8*13+1.8*\x, -3.6) { 8.\x};

\foreach \x in {19,...,22}
\node[gray] at (-1.8*18+1.8*\x, -5.4) { 8.\x};

\foreach \x in {0,...,5}
\foreach \y in {0,1,2}
\foreach \a in {0,...,7}
\draw[xshift=1.8*\x cm, yshift=-1.8*\y cm]
	(45*\a-22.5:0.8) -- (45*\a+22.5:0.8);

\foreach \x in {1,...,4}
\foreach \a in {0,...,7}
\draw[xshift=1.8*\x cm, yshift=-5.4 cm]
	(45*\a-22.5:0.8) -- (45*\a+22.5:0.8);

\foreach \a/\b/\x/\y in {
2/4/0/0, 5/7/0/0, 3/6/0/0,
2/4/1/0, 5/7/1/0,
3/6/2/0,
0/2/3/0, 4/6/3/0, 1/5/3/0,  
0/2/4/0, 4/6/4/0, 
3/5/0/1, 4/6/0/1, 2/7/0/1,  
3/5/1/1, 2/7/1/1,
2/5/2/1, 3/6/2/1, 4/7/2/1,
3/6/3/1, 
3/6/4/1, 4/7/4/1, 
3/6/5/1, 1/5/5/1,   
0/2/0/2, 1/4/0/2, 
0/2/1/2, 
1/4/2/2, 3/6/2/2,
0/4/3/2, 1/5/3/2, 2/7/3/2, 
0/4/4/2, 2/6/4/2,   
2/7/5/2,
2/7/1/3, 
2/7/2/3, 
2/7/3/3}
\draw[xshift=1.8*\x cm, yshift=-1.8*\y cm]
	(45*\a+45:0.74) to[out=45*\a+225, in=45*\b+225] (45*\b+45:0.74);

\foreach \a/\b/\x/\y in {
0/1/0/0, 
0/1/1/0, 3/6/1/0, 
0/1/2/0, 2/4/2/0, 5/7/2/0, 
3/7/3/0, 
1/5/4/0, 3/7/4/0,  
0/2/5/0, 4/6/5/0, 1/5/5/0, 3/7/5/0,
0/1/0/1, 
0/1/1/1, 4/6/1/1,
0/1/2/1, 
0/1/3/1, 2/5/3/1, 4/7/3/1, 
0/2/4/1, 1/5/4/1, 
0/2/5/1, 4/7/5/1, 
5/7/0/2, 3/6/0/2, 
5/7/1/2, 1/4/1/2, 3/6/1/2,
0/2/2/2, 5/7/2/2, 
3/6/3/2, 
1/5/4/2, 3/7/4/2,  
0/6/5/2, 1/5/5/2, 3/4/5/2, 
0/1/1/3, 4/5/1/3, 3/6/1/3,
0/1/2/3, 3/4/2/3, 5/6/2/3, 
0/5/3/3, 1/4/3/3, 3/6/3/3,  
0/7/4/3, 1/5/4/3, 2/3/4/3, 4/6/4/3}
\draw[dash pattern=on 2pt off 1pt, xshift=1.8*\x cm, yshift=-1.8*\y cm]
	(45*\a+45:0.74) to[out=45*\a+225, in=45*\b+225] (45*\b+45:0.74);
	
\foreach \a/\x/\y in {
	0/0/0, 1/0/0, 2/0/0, 5/0/0,
	0/1/0, 1/1/0, 2/1/0, 5/1/0,
	0/2/0, 1/2/0, 2/2/0, 4/2/0, 6/2/0,
	0/3/0, 3/3/0, 4/3/0, 7/3/0, 
	0/4/0, 3/4/0, 4/4/0, 7/4/0, 
	0/5/0, 2/5/0, 4/5/0, 6/5/0, 
	0/0/1, 1/0/1, 2/0/1, 
	0/1/1, 1/1/1, 2/1/1,
	0/2/1, 1/2/1, 2/2/1, 4/2/1, 6/2/1,
	0/3/1, 1/3/1, 2/3/1, 5/3/1, 
	0/4/1, 2/4/1, 4/4/1, 6/4/1, 
	0/5/1, 2/5/1, 5/5/1,   
	0/0/2, 3/0/2, 6/0/2,
	0/1/2, 3/1/2, 6/1/2, 
	0/2/2, 2/2/2, 4/2/2, 6/2/2,
	1/3/2, 3/3/2, 4/3/2, 6/3/2, 7/3/2, 
	0/4/2, 1/4/2, 4/4/2, 5/4/2,  
	1/5/2, 3/5/2, 4/5/2, 5/5/2, 7/5/2,  
	0/1/3, 1/1/3, 2/1/3,
	3/2/3, 4/2/3, 5/2/3, 6/2/3, 7/2/3,
	1/3/3, 3/3/3, 4/3/3, 6/3/3, 7/3/3, 
	2/4/3, 3/4/3, 4/4/3, 6/4/3}	
\fill[red, shift={(1.8*\x cm,-1.8*\y cm)}]
	(45*\a+22.5:0.8) circle (0.05);

\foreach \a/\x/\y in {
	3/0/0, 4/0/0, 6/0/0, 7/0/0,
	3/1/0, 4/1/0, 6/1/0, 7/1/0, 
	3/2/0, 5/2/0, 7/2/0,
	1/3/0, 2/3/0, 5/3/0, 6/3/0,
	1/4/0, 2/4/0, 5/4/0, 6/4/0, 
	1/5/0, 3/5/0, 5/5/0, 7/5/0, 
	3/0/1, 4/0/1, 5/0/1, 6/0/1, 7/0/1,
	3/1/1, 4/1/1, 5/1/1, 6/1/1, 7/1/1,
	3/2/1, 7/2/1, 5/2/1,
	3/3/1, 4/3/1, 6/3/1, 7/3/1,  
	1/4/1, 3/4/1, 5/4/1, 7/4/1, 
	1/5/1, 3/5/1, 4/5/1, 6/5/1, 7/5/1,
	1/0/2, 2/0/2, 4/0/2, 5/0/2, 7/0/2,
	1/1/2, 2/1/2, 4/1/2, 5/1/2, 7/1/2,
	1/2/2, 3/2/2, 5/2/2, 7/2/2, 
	0/3/2, 2/3/2, 5/3/2,  
	2/4/2, 3/4/2, 6/4/2, 7/4/2,  
	0/5/2, 2/5/2, 6/5/2,
	3/1/3, 4/1/3, 5/1/3, 6/1/3, 7/1/3,
	0/2/3, 1/2/3, 2/2/3, 
	0/3/3, 2/3/3, 5/3/3,  
	0/4/3, 1/4/3, 5/4/3, 7/4/3}	
\fill[blue, shift={(1.8*\x cm,-1.8*\y cm)}]
	(45*\a+22.5:0.8) circle (0.05);


\begin{scope}[yshift=-7.2cm]

\draw[gray!50]
	(-0.9,0.9) rectangle ++(5.4,-1.8)
	(-0.9,-0.9) rectangle ++(3.6,-3.6)
	(-0.9,-2.7) -- ++(3.6,0);

\foreach \x in {0,...,5}
\foreach \y in {0,...,3}
\foreach \a in {0,...,9}
\draw[xshift=1.8*\x cm, yshift=-1.8*\y cm]
	(36*\a:0.8) -- (36*\a+36:0.8);

\foreach \x in {1,...,6}
\node[gray] at (-1.8+1.8*\x, 0) { 10.\x};

\foreach \x in {7,...,12}
\node[gray] at (-1.8*7+1.8*\x, -1.8) { 10.\x};

\foreach \x in {13,...,18}
\node[gray] at (-1.8*13+1.8*\x, -3.6) { 10.\x};

\foreach \x in {19,...,24}
\node[gray] at (-1.8*19+1.8*\x, -5.4) { 10.\x};

\foreach \a/\b/\x/\y in {
	1/4/3/0, 5/8/3/0, 2/6/3/0, 3/7/3/0, 
	2/4/0/0, 5/7/0/0, 1/8/0/0, 3/6/0/0, 
	1/8/0/1, 2/5/0/1, 3/6/0/1, 4/7/0/1, 
	0/5/0/3, 1/8/0/3, 3/6/0/3, 4/9/0/3,  
	1/8/1/0, 2/4/1/0, 5/7/1/0,  
	0/4/0/2, 1/9/0/2, 3/5/0/2, 
	1/5/2/1, 4/8/2/1, 2/7/2/1,  
	6/9/3/3, 3/7/3/3, 4/8/3/3, 
	1/4/4/0, 5/8/4/0, 
	1/4/5/0, 5/8/5/0, 
	1/8/2/0, 3/6/2/0, 
	1/8/1/1, 3/6/1/1, 
	1/8/4/3, 2/5/4/3, 
	1/8/2/3, 3/6/2/3, 
	1/9/1/2, 3/5/1/2, 
	1/5/4/1, 4/8/4/1, 
	3/9/2/2, 1/5/2/2, 
	2/7/3/2, 5/9/3/2, 
	1/8/1/3, 3/6/1/3, 
	0/7/3/1, 2/9/3/1, 
	1/8/5/1, 3/9/5/1,
	1/8/5/3, 
	1/5/4/2,
	1/4/5/2 }
\draw[xshift=1.8*\x cm, yshift=-1.8*\y cm]
	(36*\a+18:0.76) to[out=36*\a+198, in=36*\b+198] (36*\b+18:0.76);

\foreach \a/\b/\x/\y in {
	0/9/3/0, 
	0/9/0/0, 
	0/9/0/1, 
	2/7/0/3, 
	0/9/1/0, 3/6/1/0, 
	6/8/0/2, 2/7/0/2, 
	0/3/2/1, 6/9/2/1,  
	0/2/3/3, 1/5/3/3, 
	0/9/4/0, 2/3/4/0, 6/7/4/0, 
	0/9/5/0, 2/7/5/0, 3/6/5/0, 
	0/9/2/0, 2/4/2/0, 5/7/2/0, 
	0/9/1/1, 2/5/1/1, 4/7/1/1, 
	0/9/4/3, 3/7/4/3, 4/6/4/3, 
	0/9/2/3, 2/7/2/3, 4/5/2/3, 
	0/4/1/2, 2/7/1/2, 6/8/1/2, 
	0/2/4/1, 7/9/4/1, 3/6/4/1, 
	0/8/2/2, 2/7/2/2, 4/6/2/2, 
	0/6/3/2, 1/3/3/2, 4/8/3/2, 
	0/4/1/3, 2/7/1/3, 5/9/1/3, 
	1/4/3/1, 3/6/3/1, 5/8/3/1, 
	0/6/5/1, 2/5/5/1, 4/7/5/1, 
	0/9/5/3, 2/6/5/3, 3/4/5/3, 5/7/5/3, 
	0/3/4/2, 2/7/4/2, 4/6/4/2, 8/9/4/2, 
	0/8/5/2, 2/7/5/2, 3/6/5/2, 5/9/5/2 }
\draw[dash pattern=on 2pt off 1pt, xshift=1.8*\x cm, yshift=-1.8*\y cm]
	(36*\a+18:0.76) to[out=36*\a+198, in=36*\b+198] (36*\b+18:0.76);
	
\foreach \a/\x/\y in {
	0/3/0, 1/3/0, 5/3/0, 9/3/0, 
	0/0/0, 1/0/0, 9/0/0, 
	0/0/1, 1/0/1, 9/0/1, 
	0/0/3, 4/0/3, 6/0/3,
	0/1/0, 1/1/0, 9/1/0,  
	2/0/2, 7/0/2, 9/0/2, 
	0/2/1, 3/2/1, 7/2/1, 
	0/3/3, 2/3/3, 6/3/3, 
	0/4/0, 1/4/0, 5/4/0, 9/4/0, 
	0/5/0, 1/5/0, 5/5/0, 9/5/0, 
	0/2/0, 1/2/0, 9/2/0, 
	0/1/1, 1/1/1, 9/1/1, 
	0/4/3, 1/4/3, 9/4/3, 
	0/2/3, 1/2/3, 9/2/3, 
	2/1/2, 7/1/2, 9/1/2, 
	0/4/1, 2/4/1, 5/4/1, 8/4/1, 
	0/2/2, 3/2/2, 8/2/2, 
	1/3/2, 3/3/2, 7/3/2,
	0/1/3, 4/1/3, 6/1/3, 
	0/3/1, 2/3/1, 5/3/1, 8/3/1, 
	0/5/1, 3/5/1, 6/5/1, 
	0/5/3, 1/5/3, 9/5/3, 
	0/4/2, 3/4/2, 8/4/2, 9/4/2, 
	0/5/2, 3/5/2, 6/5/2, 8/5/2 }	
\fill[red, shift={(1.8*\x cm,-1.8*\y cm)}]
	(36*\a:0.8) circle (0.05);

\foreach \a/\x/\y in {
	2/3/0, 4/3/0, 7/3/0,
	2/0/0, 5/0/0, 8/0/0,
	2/0/1, 4/0/1, 6/0/1, 8/0/1,
	1/0/3, 5/0/3, 9/0/3, 
	2/1/0, 5/1/0, 8/1/0,
	3/0/2, 6/0/2, 8/0/2,
	2/2/1, 5/2/1, 8/2/1,
	1/3/3, 3/3/3, 5/3/3, 8/3/3,
	2/4/0, 3/4/0, 4/4/0, 
	2/5/0, 4/5/0, 7/5/0, 
	3/2/0, 5/2/0, 7/2/0, 
	2/1/1, 5/1/1, 8/1/1,  
	3/4/3, 5/4/3, 7/4/3,
	4/2/3, 5/2/3, 6/2/3,  
	3/1/2, 6/1/2, 8/1/2, 
	1/4/1, 3/4/1, 6/4/1, 
	2/2/2, 5/2/2, 7/2/2,  
	2/3/2, 4/3/2, 8/3/2,
	1/1/3, 5/1/3, 9/1/3,     
	3/3/1, 6/3/1, 9/3/1,
	1/5/1, 4/5/1, 7/5/1, 9/5/1, 
	2/5/3, 6/5/3, 8/5/3,  
	1/4/2, 4/4/2, 6/4/2, 
	1/5/2, 5/5/2, 9/5/2 }	
\fill[blue, shift={(1.8*\x cm,-1.8*\y cm)}]
	(36*\a:0.8) circle (0.05);

\foreach \a/\x/\y in {
	3/3/0, 6/3/0, 8/3/0,
	3/0/0, 4/0/0, 6/0/0, 7/0/0,
	3/0/1, 5/0/1, 7/0/1,
	2/0/3, 3/0/3, 7/0/3, 8/0/3, 
	3/1/0, 4/1/0, 6/1/0, 7/1/0,
	0/0/2, 1/0/2, 4/0/2, 5/0/2,
	1/2/1, 4/2/1, 6/2/1, 9/2/1,
	4/3/3, 7/3/3, 9/3/3, 
	6/4/0, 7/4/0, 8/4/0,
	3/5/0, 6/5/0, 8/5/0, 
	2/2/0, 4/2/0, 6/2/0, 8/2/0,
	3/1/1, 4/1/1, 6/1/1, 7/1/1,  
	2/4/3, 4/4/3, 6/4/3, 8/4/3, 
	2/2/3, 3/2/3, 7/2/3, 8/2/3, 
	0/1/2, 1/1/2, 4/1/2, 5/1/2,
	4/4/1, 7/4/1, 9/4/1, 
	1/2/2, 4/2/2, 6/2/2, 9/2/2, 
	0/3/2, 5/3/2, 6/3/2, 9/3/2,
	2/1/3, 3/1/3, 7/1/3, 8/1/3, 
	1/3/1, 4/3/1, 7/3/1,
	2/5/1, 5/5/1, 8/5/1,  
	3/5/3, 4/5/3, 5/5/3, 7/5/3,
	2/4/2, 5/4/2, 7/4/2, 
	2/5/2, 4/5/2, 7/5/2 }	
\fill[green, shift={(1.8*\x cm,-1.8*\y cm)}]
	(36*\a:0.8) circle (0.05);
	
\end{scope}


\begin{scope}[yshift=-14.4cm]

\foreach \a in {0,...,11}
{
\foreach \x in {0,...,4}
\draw[xshift=0.9cm + 1.8*\x cm]
	(30*\a-15:0.8) -- (30*\a+15:0.8);	

\foreach \x in {0,...,5}
\draw[shift={(0.9+1.8*\x cm, -1.8 cm)}]
	(30*\a-15:0.8) -- (30*\a+15:0.8);
}

\foreach \x in {1,...,5}
\node[gray] at (-0.9+1.8*\x, 0) { 12.\x};

\foreach \x in {6,...,11}
\node[gray] at (-1.8*6+1.8*\x, -1.8) { 12.\x};
   
\foreach \a/\b/\x/\y in {
	0/5/0/0, 1/10/0/0, 3/6/0/0, 4/11/0/0, 
	1/10/1/0, 3/6/1/0,
	1/10/2/0, 3/6/2/0, 
	3/11/3/0, 1/5/3/0,
	3/10/4/0, 1/6/4/0,
	0/4/0/1, 1/9/0/1, 3/10/0/1, 
	2/11/1/1, 3/6/1/1, 7/10/1/1, 4/8/1/1, 5/9/1/1,
	2/11/2/1, 3/6/2/1, 7/10/2/1,
	2/11/3/1, 3/6/3/1, 7/10/3/1,
	0/5/4/1, 1/8/4/1, 4/9/4/1, 
	2/11/5/1, 4/8/5/1
	}
\draw[shift={(0.9 cm - 0.9*\y cm + 1.8*\x cm, -1.8*\y cm)}]
	(30*\a+30:0.77) to[out=30*\a+210, in=30*\b+210] (30*\b+30:0.77);

\foreach \a/\b/\x/\y in {
	2/8/0/0, 7/9/0/0, 
	7/9/1/0, 0/4/1/0, 2/8/1/0, 5/11/1/0, 
	0/11/2/0, 2/8/2/0, 4/5/2/0, 7/9/2/0, 
	0/9/3/0, 2/8/3/0, 4/7/3/0, 6/10/3/0, 
	0/4/4/0, 2/8/4/0, 5/7/4/0, 9/11/4/0, 
	2/6/0/1, 5/8/0/1, 7/11/0/1,   
	0/1/1/1, 
	0/1/2/1, 5/8/2/1, 4/9/2/1,
	0/1/3/1, 4/5/3/1, 8/9/3/1, 
	2/11/4/1, 3/6/4/1, 7/10/4/1, 
	0/1/5/1, 7/9/5/1, 3/6/5/1, 5/10/5/1
	}
\draw[dash pattern=on 2pt off 1pt, shift={(0.9 cm - 0.9*\y cm + 1.8*\x cm, -1.8*\y cm)}]
	(30*\a+30:0.77) to[out=30*\a+210, in=30*\b+210] (30*\b+30:0.77);

\foreach \a/\x/\y in {
	3/0/0, 7/0/0, 9/0/0, 
	3/1/0, 7/1/0, 9/1/0, 
	3/2/0, 7/2/0, 9/2/0, 
	0/3/0, 3/3/0, 9/3/0, 
	3/4/0, 9/4/0, 11/4/0, 
	1/0/1, 4/0/1, 10/0/1,
	3/1/1, 7/1/1, 11/1/1,
	3/2/1, 7/2/1, 11/2/1,
	3/3/1, 7/3/1, 11/3/1,
	1/4/1, 5/4/1, 9/4/1,
	0/5/1, 1/5/1, 2/5/1
	}
\fill[shift={(0.9 cm - 0.9*\y cm + 1.8*\x cm,-1.8*\y cm)}]
	(30*\a+15:0.8) circle (0.05);

\foreach \a/\x/\y in {
	2/0/0, 8/0/0, 10/0/0, 
	2/1/0, 8/1/0, 10/1/0, 
	2/2/0, 8/2/0, 10/2/0, 
	2/3/0, 5/3/0, 8/3/0,
	2/4/0, 6/4/0, 8/4/0, 
	2/0/1, 6/0/1, 9/0/1,
	0/1/1, 1/1/1, 2/1/1,
	0/2/1, 1/2/1, 2/2/1,
	0/3/1, 1/3/1, 2/3/1,
	2/4/1, 8/4/1, 11/4/1,
	3/5/1, 6/5/1, 11/5/1
	}
\fill[red, shift={(0.9 cm - 0.9*\y cm + 1.8*\x cm,-1.8*\y cm)}]
	(30*\a+15:0.8) circle (0.05);

\foreach \a/\x/\y in {  
	0/0/0, 4/0/0, 6/0/0,
	0/1/0, 4/1/0, 6/1/0, 
	4/2/0, 5/2/0, 6/2/0,
	11/3/0, 4/3/0, 7/3/0,
	1/4/0, 5/4/0, 7/4/0,
	0/0/1, 5/0/1, 8/0/1,
	4/1/1, 6/1/1, 9/1/1,
	4/2/1, 6/2/1, 9/2/1,
	4/3/1, 5/3/1, 6/3/1, 
	0/4/1, 3/4/1, 6/4/1,
	4/5/1, 7/5/1, 9/5/1
	}
\fill[blue, shift={(0.9 cm - 0.9*\y cm + 1.8*\x cm,-1.8*\y cm)}]
	(30*\a+15:0.8) circle (0.05);
	
\foreach \a/\x/\y in {
	1/0/0, 5/0/0, 11/0/0, 
	1/1/0, 5/1/0, 11/1/0,
	0/2/0, 1/2/0, 11/2/0, 
	1/3/0, 6/3/0, 10/3/0,
	0/4/0, 4/4/0, 10/4/0,
	3/0/1, 7/0/1, 11/0/1,
	5/1/1, 8/1/1, 10/1/1,
	5/2/1, 8/2/1, 10/2/1,
	8/3/1, 9/3/1, 10/3/1,
	4/4/1, 7/4/1, 10/4/1,
	5/5/1, 8/5/1, 10/5/1
	}
\fill[green, shift={(0.9 cm - 0.9*\y cm + 1.8*\x cm,-1.8*\y cm)}]
	(30*\a+15:0.8) circle (0.05);

\end{scope}
		
\end{tikzpicture}
\caption{Tilings of $3{\bb P}^2$ by single $n$-gon, $n\ge 7$.} 
\label{3P2_n8}
\end{figure}

\section{Geometrical Realisation}

The discussion so far is combinatorial. The geometrical realisation means the tiling is by a single polygon with straight line edges, Specifically, for any planar diagram $D$, we need to find a polygon in the Poincar\'e disk ${\bb H}^2$, such that all edges are straight (i.e., geodesic), and the paired edges have the same length, and the sum of the angle values of the corners at any vertex is $2\pi$. 

We denote the {\em angle value} of the $i$-th corner $i$ by $[i]$, and denote the {\em length} of the $i$-th edge $\bar{i}$ by $|\bar{i}|$. Then the geometric condition on the decagon in Figures \ref{edge_pair} and \ref{pd_vs} is 
\begin{align*}
|\bar{0}| &=|\bar{2}|,\;
|\bar{1}|=|\bar{4}|,\;
|\bar{3}|=|\bar{7}|,\;
|\bar{5}|=|\bar{8}|,\;
|\bar{6}|=|\bar{9}|, \\
[0] &+[3]+[6]+[8]=2\pi,\;
[1]+[2]+[4]+[5]+[7]+[9]=2\pi.
\end{align*}
The geometric condition on the decagon in Figure \ref{edge_pair3} is
\begin{align*}
|\bar{0}| &=|\bar{1}|,\;
|\bar{2}|=|\bar{5}|,\;
|\bar{3}|=|\bar{8}|,\;
|\bar{4}|=|\bar{7}|,\;
|\bar{6}|=|\bar{9}|, \\
[0] &+[1]+[2]+[6]=2\pi,\;
[3]+[5]+[8]=2\pi,\;
[4]+[7]+[9]=2\pi.
\end{align*}

The geometric condition on the tilings 8.1 and 8.2 in Figure \ref{3P2_n8} is
\begin{align}
|\bar{0}| &=|\bar{1}|,\;
|\bar{2}|=|\bar{4}|,\;
|\bar{3}|=|\bar{6}|,\;
|\bar{5}|=|\bar{7}|, \label{oct1E} \\
[0] &+[1]+[2]+[5]=[3]+[4]+[6]+[7]=2\pi. \label{oct1A}
\end{align}
Since the two planar diagrams have the same geometric conditions, an octagon tiles $3{\bb P}^2$ in the 8.1 way if and only if it tiles $3{\bb P}^2$ in the 8.2 way. In other words, such an octagon can tile $3{\bb P}^2$ in exactly two ways. The same happens to the following pairs of tilings in Figure \ref{3P2_n8}: 8.4 and 8.5, 8.7 and 8.8, 8.13 and 8.14, 10.1 and 10.2, 10.13 and 10.14. 

On the other hand, if the surface is orientable, then all pairs are opposing. Therefore different tilings have different pairings. In particular, a general polygon for one tiling should not be suitable for another tiling. 

We believe all the combinatorial single tile tilings have geometrical realisations. However, we do not have a general proof. Of course, if all the vertices have the same degree, then regular polygon provides geometrical realisation. This includes the two extreme cases at the end of Section \ref{intro}:
\begin{itemize}
\item $n=2(2-\chi)$ is the minimal value, which means the tiling has a single vertex. An example is the tiling of $2{\bb T}^2$ by single octagon in Figure \ref{2T2_n8}.
\item $n=6(1-\chi)$ is the maximal value, which means all vertices have degree $3$.
\end{itemize}
This also includes ten tilings in Figure \ref{3P2_n8} (8.1, 8.2, 8.4, 8.5, 8.6, etc.) in which all vertices have degree 4. 

In general, we may show the geometrical realisations by actually finding specific examples satisfying the edge length and angle sum equalities. Of course the method can only be done one planar diagram at a time. Moreover, we wish to show that a planar diagram has geometrical realisation that does not satisfy different geometrical conditions from different planar diagrams. Therefore it is still worthwhile to do such case by case calculations even when regular polygons already provide geometrical realisations.

The governing equation for a hyperbolic $n$-gon is
\begin{equation}\label{hpoly}
\prod_{i=0}^{n-1}R(\pi-[i])S(|\bar{i}|)=I.
\end{equation}
The equation appeared in Section 3.3 of \cite{cly1} and Section 3.4 of \cite{cly2}. We illustrate the equation by the quadrilateral in Figure \ref{polygon}. The gray dot is the center of the Poincar\'e disk ${\bb H}^2$, and the arrow is horizontal and to the right. In the first picture, we put the quadrilateral in the position such that the corner 0 is the center of the Poincar\'e disk, and the edge $\bar{3}$ is horizontal and on the left of the center. Next, applying the horizontal shift $S(|\bar{3}|)$ that moves the quadrilateral to the right by the length $|\bar{3}|$, we get the second picture, where the corner 3 is at the center, and the edge $\bar{3}$ is still horizontal. Then applying the rotation $R(\pi-[3])$ around the center, we get the third picture, where the corner 3 is at the center, and the edge $\bar{2}$ is horizontal and on the left of the center. Then applying the horizontal shift $S(|\bar{2}|)$, and so on. All the moves and shifts that correspond to the edges and corners of the quadrilateral combine to form the transformation $R(\pi-[0])S(|\bar{0}|)R(\pi-[1])S(|\bar{1}|)R(\pi-[2])S(|\bar{2}|)R(\pi-[3])S(|\bar{3}|)$. Since applying the transformation to the quadrilateral gives back the quadrilateral at the beginning, the composition is the identity.

\begin{figure}[htp]
\centering
\begin{tikzpicture}[>=latex]

\foreach \x in {0,...,3}
{
\draw[->, gray!50]
	(-1.4+3*\x,0) -- ++(2.8,0);

\fill[gray!50] (3*\x,0) circle (0.07);
}

\draw[xshift=-0.9 cm]
	(0,0) -- ++(130:1) -- ++(30:1.2) -- (0:0.9) -- cycle;

\draw[xshift=3 cm]
	(0,0) -- ++(130:1) -- ++(30:1.2) -- (0:0.9) -- cycle;
		
\draw[xshift=6cm, rotate=50]
	(0,0) -- ++(130:1) -- ++(30:1.2) -- (0:0.9) -- cycle;

\draw[xshift=6cm, rotate=50]
	(0,0) -- ++(130:1) -- ++(30:1.2) -- (0:0.9) -- cycle;

\draw[xshift=10cm, rotate=50]
	(0,0) -- ++(130:1) -- ++(30:1.2) -- (0:0.9) -- cycle;

\draw[xshift=3 cm]
	(-0.3,0) arc (180:130:0.3);
\draw[xshift=9 cm]
	(-0.3,0) arc (180:80:0.3);
	
\foreach \x/\y/\n in {
	-0.25/0.2/0, -0.6/1.1/1, -1.25/0.73/2, -0.8/0.2/3, 
	3.65/0.2/0, 3.3/1.1/1, 2.65/0.73/2, 3.1/0.2/3, 
	6.25/0.6/0, 5.35/0.9/1, 5.25/0.2/2, 5.9/0.2/3, 
	10.25/0.6/0, 9.35/0.9/1, 9.25/0.2/2, 9.9/0.2/3
	}	
\node at (\x,\y) {\scriptsize \n};

\node at (-0.5,-0.15) {\tiny $|\bar{3}|$};
\node at (2.35,0.2) {\tiny $\pi\!-\![3]$};
\node at (5.5,-0.15) {\tiny $|\bar{2}|$};
\node at (8.4,0.25) {\tiny $\pi\!-\![2]$};

\draw[->, very thick]
	(0.7,0.7) -- ++(1,0);
\node at (1.2,1) {\scriptsize $S(|\bar{3}|)$};

\draw[->, very thick]
	(3.9,0.7) -- ++(1,0);
\node at (4.4,1) {\scriptsize $R(\pi\!-\![3])$};

\draw[->, very thick]
	(7.4,0.7) -- ++(1,0);
\node at (7.9,1) {\scriptsize $S(|\bar{2}|)$};
		
\end{tikzpicture}
\caption{A polygon is a sequence of rotations.}
\label{polygon}
\end{figure}

The isometry group of the Poincar\'e disk is $PSU(1,1)$. In terms of the complex numbers in the unit disk, the two transformations are
\[
S(l)\colon z\mapsto \frac{z\cosh\frac{l}{2}+\sinh\frac{l}{2}}{z\sinh\frac{l}{2}+\cosh\frac{l}{2}};\quad
R(\theta)\colon z\mapsto e^{i\theta}z.
\]
It is easier to use the faithful $3$-dimensional irreducible representation of $PSU(1,1)$, for which the two transformations are represented by the matrices
\[
S(l)=\begin{pmatrix}
\cosh l & 0 & \sinh l \\
0 & 1 & 0 \\
\sinh l & 0 & \cosh l
\end{pmatrix},\quad
R(\theta)=\begin{pmatrix}
\cos \theta & -\sin \theta & 0 \\
\sin \theta & \cos \theta & 0 \\
0 & 0 & 1
\end{pmatrix}.
\]

For an $n$-gon, the number of variables is $\frac{n}{2}$ edge lengths and $n$ angle values. The total number is $\frac{3}{2}n$. On the other hand, the equations are the angle sums at vertices and the equality \eqref{hpoly} in $PSU(1,1)$. The number of vertices is $v=\chi+e-f=\chi+\frac{n}{2}-1$. The dimension is $PSU(1,1)$ is 3. The total number is $v+3=\chi+\frac{n}{2}+2$. Therefore the number of free variables in the $n$-gon is $\frac{3}{2}n-(\chi+\frac{n}{2}+2)=n-2-\chi$. Using the bound $2(2-\chi)\le n\le 6(1-\chi)$, we know the number of variables is between $2-3\chi$ and $4-7\chi$.

The equation \eqref{hpoly} essentially means that line segments of certain length and the turning angles determine a closed loop. If the closed loop is simple, i.e., no self intersection, then the line segments and turning angles determine a {\em simple} hyperbolic $n$-gon. The polygons in tilings are required to be simple. To check that a solution to \eqref{hpoly} gives a simple polygon, we draw the polygon and visually verify the simple property.

The geometric condition has the edge part and the angle part. For tilings 8.1 and 8.2 of Figure \ref{3P2_n8}, we first set the edge part \eqref{oct1E} as
\[
|\bar{0}|=|\bar{1}|=10a,\;
|\bar{2}|=|\bar{4}|=11a,\;
|\bar{3}|=|\bar{6}|=12a,\;
|\bar{5}|=|\bar{7}|=13a.
\]
The idea is to set different edge pairs to have different lengths, such that the polygon is not suitable for any tiling with different edge pairing. In particular, a polygon satisfying the edge equalities above is suitable only for tilings 8.1, 8.2, 8.3 in Figure \ref{3P2_n8}. Moreover, we use the multiples 10, 11, 12, 13 instead of 1, 2, 3, 4 because we expect it to be harder to get simple polygon if the edge lengths vary too much. 

Next we choose angle values carefully, such that the polygon is suitable for tilings 8.1 and 8.2 of Figure \ref{3P2_n8}, and is not suitable for tiling 8.3 of Figure \ref{3P2_n8}. This can be done if the angle $[5]$ in  8.1 and 8.2 is not equal to $[4]+[6]$ in 8.3. For example, we may set the following for 8.1 and 8.2
\[
[3]=[4]=[6]=[7]=[5]=\tfrac{\pi}{2},\;
[1]=x,\;
[2]=y,\;
[0]=\tfrac{3\pi}{2}-x-y, 
\]
and set the following for 8.3
\[
[3]=[7]=[5]=\tfrac{2\pi}{3},\;
[4]=[6]=\tfrac{\pi}{5},\;
[1]=x,\;
[2]=y,\;
[0]=\tfrac{8\pi}{5}-x-y. 
\]
The reason we allow three variables $a,x,y$ is that the dimension of $PSU(1,1)$ is 3.

We substitute the settings into the equation \eqref{hpoly}, and calculate the polygon. The solution for 8.1 and 8.2 is
\[
a=0.14332,\;
[1]=0.88279\pi,\;
[2]=0.35894\pi,
\]
and the solution for 8.3 is 
\[
a=0.19302,\;
[1]=1.37208\pi,\;
[2]=0.15045\pi.
\]
Then we draw the octagons, which are the first three in Figure \ref{3P2_n8geom}. The point $\circ$ indicates the vertex of the octagon on the lower left of Figure \ref{3P2_n8geom}, compared with the planar diagrams $8.*$ in Figure \ref{3P2_n8}. 

We carry out the similar calculations for the other tilings in Figure \ref{3P2_n8}, and get the corresponding octagons in Figure \ref{3P2_n8geom}. The parameters are all set to make sure that the octagon for one tiling cannot be used for tilings with different edge length and angle equalities. Since all octagons are simple, we obtain the geometric realisations.

\begin{figure}[htp]
\centering
\begin{tikzpicture}[scale=1] 

\foreach \x in {1,...,6}
\node[gray] at (-2.3+2.3*\x, 0) { 8.\x};

\foreach \x in {7,...,12}
\node[gray] at (-2.3*7+2.3*\x, -2.3) { 8.\x};

\foreach \x in {13,...,18}
\node[gray] at (-2.3*13+2.3*\x, -4.6) { 8.\x};

\foreach \x in {19,...,22}
\node[gray] at (-2.3*18+2.3*\x, -6.9) { 8.\x};

\foreach \x in {0,...,5}
\foreach \y in {0,1,2}
\draw[shift={(2.3*\x,-2.3*\y)}] (0,0) circle (1);

\foreach \x in {1,2,3,4}
\draw[shift={(2.3*\x,-2.3*3)}] (0,0) circle (1);
\begin{scope}[shift={(0,0)}]

\drawhypgeodesic{-0.01029}{-0.78132}{0.26865}{-0.30340}
\drawhypgeodesic{0.26865}{-0.30340}{0.68187}{0.17257}
\drawhypgeodesic{0.68187}{0.17257}{0.44858}{0.54489}
\drawhypgeodesic{0.44858}{0.54489}{0.05957}{0.67003}
\drawhypgeodesic{0.05957}{0.67003}{-0.40438}{0.44211}
\drawhypgeodesic{-0.40438}{0.44211}{-0.61457}{-0.18241}
\drawhypgeodesic{-0.61457}{-0.18241}{-0.42943}{-0.56246}
\drawhypgeodesic{-0.42943}{-0.56246}{-0.01029}{-0.78132}

\filldraw[fill=white] (-0.01029,-0.78132) circle (0.05);  
\end{scope}
\begin{scope}[shift={(2.3,0)}]

\drawhypgeodesic{-0.01029}{-0.78132}{0.26865}{-0.30340}
\drawhypgeodesic{0.26865}{-0.30340}{0.68187}{0.17257}
\drawhypgeodesic{0.68187}{0.17257}{0.44858}{0.54489}
\drawhypgeodesic{0.44858}{0.54489}{0.05957}{0.67003}
\drawhypgeodesic{0.05957}{0.67003}{-0.40438}{0.44211}
\drawhypgeodesic{-0.40438}{0.44211}{-0.61457}{-0.18241}
\drawhypgeodesic{-0.61457}{-0.18241}{-0.42943}{-0.56246}
\drawhypgeodesic{-0.42943}{-0.56246}{-0.01029}{-0.78132}

\filldraw[fill=white] (-0.01029,-0.78132) circle (0.05);  
\end{scope}
\begin{scope}[shift={(4.6,0)}]

\drawhypgeodesic{0.17009}{-0.86087}{0.04632}{-0.19868}
\drawhypgeodesic{0.04632}{-0.19868}{0.79964}{0.13228}
\drawhypgeodesic{0.79964}{0.13228}{0.38683}{0.52110}
\drawhypgeodesic{0.38683}{0.52110}{0.04322}{0.83646}
\drawhypgeodesic{0.04322}{0.83646}{-0.34950}{0.39542}
\drawhypgeodesic{-0.34950}{0.39542}{-0.76330}{-0.28295}
\drawhypgeodesic{-0.76330}{-0.28295}{-0.33330}{-0.54277}
\drawhypgeodesic{-0.33330}{-0.54277}{0.17009}{-0.86087}

\filldraw[fill=white] (0.17009,-0.86087) circle (0.05);      
\end{scope}
\begin{scope}[shift={(6.9,0)}]

\drawhypgeodesic{-0.14220}{-0.56128}{0.35594}{-0.36590}
\drawhypgeodesic{0.35594}{-0.36590}{0.75978}{-0.13007}
\drawhypgeodesic{0.75978}{-0.13007}{0.65822}{0.20689}
\drawhypgeodesic{0.65822}{0.20689}{0.21053}{0.54135}
\drawhypgeodesic{0.21053}{0.54135}{-0.44796}{0.43247}
\drawhypgeodesic{-0.44796}{0.43247}{-0.76532}{0.16378}
\drawhypgeodesic{-0.76532}{0.16378}{-0.62899}{-0.28725}
\drawhypgeodesic{-0.62899}{-0.28725}{-0.14220}{-0.56128}

\filldraw[fill=white] (-0.14220,-0.56128) circle (0.05);      
\end{scope}
\begin{scope}[shift={(9.2,0)}]

\drawhypgeodesic{-0.14220}{-0.56128}{0.35594}{-0.36590}
\drawhypgeodesic{0.35594}{-0.36590}{0.75978}{-0.13007}
\drawhypgeodesic{0.75978}{-0.13007}{0.65822}{0.20689}
\drawhypgeodesic{0.65822}{0.20689}{0.21053}{0.54135}
\drawhypgeodesic{0.21053}{0.54135}{-0.44796}{0.43247}
\drawhypgeodesic{-0.44796}{0.43247}{-0.76532}{0.16378}
\drawhypgeodesic{-0.76532}{0.16378}{-0.62899}{-0.28725}
\drawhypgeodesic{-0.62899}{-0.28725}{-0.14220}{-0.56128}

\filldraw[fill=white] (-0.14220,-0.56128) circle (0.05);      
\end{scope}
\begin{scope}[shift={(11.5,0)}]
\drawhypgeodesic{-0.30537}{-0.57115}{0.19331}{-0.39336}
\drawhypgeodesic{0.19331}{-0.39336}{0.66223}{-0.04977}
\drawhypgeodesic{0.66223}{-0.04977}{0.74470}{0.28824}
\drawhypgeodesic{0.74470}{0.28824}{0.38650}{0.53478}
\drawhypgeodesic{0.38650}{0.53478}{-0.28630}{0.44055}
\drawhypgeodesic{-0.28630}{0.44055}{-0.69267}{0.10843}
\drawhypgeodesic{-0.69267}{0.10843}{-0.70240}{-0.35772}
\drawhypgeodesic{-0.70240}{-0.35772}{-0.30537}{-0.57115}

\filldraw[fill=white] (-0.30537,-0.57115) circle (0.05);      
\end{scope}
\begin{scope}[shift={(0,-2.3)}]
\drawhypgeodesic{-0.11187}{-0.71457}{0.01710}{-0.10403}
\drawhypgeodesic{0.01710}{-0.10403}{0.64454}{-0.26552}
\drawhypgeodesic{0.64454}{-0.26552}{0.69662}{0.22977}
\drawhypgeodesic{0.69662}{0.22977}{0.31762}{0.62922}
\drawhypgeodesic{0.31762}{0.62922}{-0.33685}{0.59785}
\drawhypgeodesic{-0.33685}{0.59785}{-0.68749}{0.11654}
\drawhypgeodesic{-0.68749}{0.11654}{-0.53968}{-0.48926}
\drawhypgeodesic{-0.53968}{-0.48926}{-0.11187}{-0.71457}

\filldraw[fill=white] (-0.11187,-0.71457) circle (0.05);      
\end{scope}
\begin{scope}[shift={(2.3,-2.3)}]
\drawhypgeodesic{-0.11187}{-0.71457}{0.01710}{-0.10403}
\drawhypgeodesic{0.01710}{-0.10403}{0.64454}{-0.26552}
\drawhypgeodesic{0.64454}{-0.26552}{0.69662}{0.22977}
\drawhypgeodesic{0.69662}{0.22977}{0.31762}{0.62922}
\drawhypgeodesic{0.31762}{0.62922}{-0.33685}{0.59785}
\drawhypgeodesic{-0.33685}{0.59785}{-0.68749}{0.11654}
\drawhypgeodesic{-0.68749}{0.11654}{-0.53968}{-0.48926}
\drawhypgeodesic{-0.53968}{-0.48926}{-0.11187}{-0.71457}

\filldraw[fill=white] (-0.11187,-0.71457) circle (0.05);      
\end{scope}
\begin{scope}[shift={(4.6,-2.3)}]
\drawhypgeodesic{0.34758}{-0.58689}{0.65360}{-0.61259}
\drawhypgeodesic{0.65360}{-0.61259}{0.56190}{-0.30254}
\drawhypgeodesic{0.56190}{-0.30254}{0.30359}{0.36919}
\drawhypgeodesic{0.30359}{0.36919}{-0.09880}{0.81357}
\drawhypgeodesic{-0.09880}{0.81357}{-0.51266}{0.42772}
\drawhypgeodesic{-0.51266}{0.42772}{-0.82965}{0.13980}
\drawhypgeodesic{-0.82965}{0.13980}{-0.42557}{-0.24825}
\drawhypgeodesic{-0.42557}{-0.24825}{0.34758}{-0.58689}

\filldraw[fill=white] (0.34758,-0.58689) circle (0.05);      
\end{scope}
\begin{scope}[shift={(6.9,-2.3)}]
\drawhypgeodesic{0.11370}{-0.74394}{0.24776}{-0.24553}
\drawhypgeodesic{0.24776}{-0.24553}{0.74228}{-0.08017}
\drawhypgeodesic{0.74228}{-0.08017}{0.55231}{0.35158}
\drawhypgeodesic{0.55231}{0.35158}{0.00052}{0.62580}
\drawhypgeodesic{0.00052}{0.62580}{-0.57185}{0.50320}
\drawhypgeodesic{-0.57185}{0.50320}{-0.67071}{0.06519}
\drawhypgeodesic{-0.67071}{0.06519}{-0.41399}{-0.47612}
\drawhypgeodesic{-0.41399}{-0.47612}{0.11370}{-0.74394}

\filldraw[fill=white] (0.11370,-0.74394) circle (0.05);      
\end{scope}
\begin{scope}[shift={(9.2,-2.3)}]
\drawhypgeodesic{0.10877}{-0.64247}{0.47572}{-0.46600}
\drawhypgeodesic{0.47572}{-0.46600}{0.67634}{-0.07522}
\drawhypgeodesic{0.67634}{-0.07522}{0.54914}{0.31303}
\drawhypgeodesic{0.54914}{0.31303}{0.07802}{0.59827}
\drawhypgeodesic{0.07802}{0.59827}{-0.48776}{0.45635}
\drawhypgeodesic{-0.48776}{0.45635}{-0.71517}{0.11209}
\drawhypgeodesic{-0.71517}{0.11209}{-0.68508}{-0.29605}
\drawhypgeodesic{-0.68508}{-0.29605}{0.10877}{-0.64247}

\filldraw[fill=white] (0.10877,-0.64247) circle (0.05);      
\end{scope}
\begin{scope}[shift={(11.5,-2.3)}]
\drawhypgeodesic{-0.10667}{-0.54712}{0.34359}{-0.59779}
\drawhypgeodesic{0.34359}{-0.59779}{0.56790}{-0.15896}
\drawhypgeodesic{0.56790}{-0.15896}{0.64080}{0.27183}
\drawhypgeodesic{0.64080}{0.27183}{0.29174}{0.61708}
\drawhypgeodesic{0.29174}{0.61708}{-0.31640}{0.44979}
\drawhypgeodesic{-0.31640}{0.44979}{-0.71425}{0.19331}
\drawhypgeodesic{-0.71425}{0.19331}{-0.70671}{-0.22814}
\drawhypgeodesic{-0.70671}{-0.22814}{-0.10667}{-0.54712}

\filldraw[fill=white] (-0.10667,-0.54712) circle (0.05);      
\end{scope}
\begin{scope}[shift={(0,-4.6)}]
\drawhypgeodesic{-0.03316}{-0.58635}{0.37844}{-0.64366}
\drawhypgeodesic{0.37844}{-0.64366}{0.59724}{-0.26867}
\drawhypgeodesic{0.59724}{-0.26867}{0.44902}{0.21502}
\drawhypgeodesic{0.44902}{0.21502}{0.05799}{0.68057}
\drawhypgeodesic{0.05799}{0.68057}{-0.32582}{0.72081}
\drawhypgeodesic{-0.32582}{0.72081}{-0.52147}{0.22349}
\drawhypgeodesic{-0.52147}{0.22349}{-0.60225}{-0.34119}
\drawhypgeodesic{-0.60225}{-0.34119}{-0.03316}{-0.58635}

\filldraw[fill=white] (-0.03316,-0.58635) circle (0.05);      
\end{scope}
\begin{scope}[shift={(2.3,-4.6)}]
\drawhypgeodesic{-0.03316}{-0.58635}{0.37844}{-0.64366}
\drawhypgeodesic{0.37844}{-0.64366}{0.59724}{-0.26867}
\drawhypgeodesic{0.59724}{-0.26867}{0.44902}{0.21502}
\drawhypgeodesic{0.44902}{0.21502}{0.05799}{0.68057}
\drawhypgeodesic{0.05799}{0.68057}{-0.32582}{0.72081}
\drawhypgeodesic{-0.32582}{0.72081}{-0.52147}{0.22349}
\drawhypgeodesic{-0.52147}{0.22349}{-0.60225}{-0.34119}
\drawhypgeodesic{-0.60225}{-0.34119}{-0.03316}{-0.58635}

\filldraw[fill=white] (-0.03316,-0.58635) circle (0.05);      
\end{scope}
\begin{scope}[shift={(4.6,-4.6)}]
\drawhypgeodesic{-0.11714}{-0.62773}{0.31373}{-0.47487}
\drawhypgeodesic{0.31373}{-0.47487}{0.63805}{-0.08405}
\drawhypgeodesic{0.63805}{-0.08405}{0.57193}{0.32675}
\drawhypgeodesic{0.57193}{0.32675}{0.15295}{0.62610}
\drawhypgeodesic{0.15295}{0.62610}{-0.31721}{0.55147}
\drawhypgeodesic{-0.31721}{0.55147}{-0.64715}{0.07728}
\drawhypgeodesic{-0.64715}{0.07728}{-0.59517}{-0.39496}
\drawhypgeodesic{-0.59517}{-0.39496}{-0.11714}{-0.62773}

\filldraw[fill=white] (-0.11714,-0.62773) circle (0.05);      
\end{scope}
\begin{scope}[shift={(6.9,-4.6)}]
\drawhypgeodesic{0.00959}{-0.62607}{0.41704}{-0.57515}
\drawhypgeodesic{0.41704}{-0.57515}{0.56941}{-0.12136}
\drawhypgeodesic{0.56941}{-0.12136}{0.54839}{0.42300}
\drawhypgeodesic{0.54839}{0.42300}{0.08381}{0.68671}
\drawhypgeodesic{0.08381}{0.68671}{-0.30125}{0.44076}
\drawhypgeodesic{-0.30125}{0.44076}{-0.69495}{0.13265}
\drawhypgeodesic{-0.69495}{0.13265}{-0.63203}{-0.36054}
\drawhypgeodesic{-0.63203}{-0.36054}{0.00959}{-0.62607}

\filldraw[fill=white] (0.00959,-0.62607) circle (0.05);      
\end{scope}
\begin{scope}[shift={(9.2,-4.6)}]
\drawhypgeodesic{0.15353}{-0.56971}{0.53118}{-0.59378}
\drawhypgeodesic{0.53118}{-0.59378}{0.66607}{-0.25095}
\drawhypgeodesic{0.66607}{-0.25095}{0.49343}{0.27968}
\drawhypgeodesic{0.49343}{0.27968}{-0.15353}{0.56971}
\drawhypgeodesic{-0.15353}{0.56971}{-0.53118}{0.59378}
\drawhypgeodesic{-0.53118}{0.59378}{-0.66607}{0.25095}
\drawhypgeodesic{-0.66607}{0.25095}{-0.49343}{-0.27968}
\drawhypgeodesic{-0.49343}{-0.27968}{0.15353}{-0.56971}

\filldraw[fill=white] (0.15353,-0.56971) circle (0.05);      
\end{scope}
\begin{scope}[shift={(11.5,-4.6)}]
\drawhypgeodesic{0.01930}{-0.66679}{0.40747}{-0.57382}
\drawhypgeodesic{0.40747}{-0.57382}{0.50402}{-0.06707}
\drawhypgeodesic{0.50402}{-0.06707}{0.46218}{0.51412}
\drawhypgeodesic{0.46218}{0.51412}{0.01058}{0.73815}
\drawhypgeodesic{0.01058}{0.73815}{-0.44693}{0.53171}
\drawhypgeodesic{-0.44693}{0.53171}{-0.49177}{0.00539}
\drawhypgeodesic{-0.49177}{0.00539}{-0.46485}{-0.48169}
\drawhypgeodesic{-0.46485}{-0.48169}{0.01930}{-0.66679}

\filldraw[fill=white] (0.01930,-0.66679) circle (0.05);      
\end{scope}
\begin{scope}[shift={(2.3,-6.9)}]
\drawhypgeodesic{0.62629}{-0.55817}{0.79931}{-0.42000}
\drawhypgeodesic{0.79931}{-0.42000}{0.47952}{-0.37117}
\drawhypgeodesic{0.47952}{-0.37117}{0.52779}{0.30063}
\drawhypgeodesic{0.52779}{0.30063}{-0.06699}{0.67057}
\drawhypgeodesic{-0.06699}{0.67057}{-0.63762}{0.45613}
\drawhypgeodesic{-0.63762}{0.45613}{-0.85598}{0.09887}
\drawhypgeodesic{-0.85598}{0.09887}{-0.87216}{-0.17689}
\drawhypgeodesic{-0.87216}{-0.17689}{0.62629}{-0.55817}

\filldraw[fill=white] (0.62629,-0.55817) circle (0.05);      
\end{scope}
\begin{scope}[shift={(4.6,-6.9)}]
\drawhypgeodesic{-0.09952}{-0.74084}{-0.03415}{-0.14678}
\drawhypgeodesic{-0.03415}{-0.14678}{0.62746}{-0.23603}
\drawhypgeodesic{0.62746}{-0.23603}{0.67298}{0.28151}
\drawhypgeodesic{0.67298}{0.28151}{0.29252}{0.66330}
\drawhypgeodesic{0.29252}{0.66330}{-0.27747}{0.63096}
\drawhypgeodesic{-0.27747}{0.63096}{-0.67248}{0.07982}
\drawhypgeodesic{-0.67248}{0.07982}{-0.50934}{-0.53194}
\drawhypgeodesic{-0.50934}{-0.53194}{-0.09952}{-0.74084}

\filldraw[fill=white] (-0.09952,-0.74084) circle (0.05);      
\end{scope}
\begin{scope}[shift={(6.9,-6.9)}]
\drawhypgeodesic{-0.25001}{-0.51389}{0.28106}{-0.57445}
\drawhypgeodesic{0.28106}{-0.57445}{0.61389}{-0.03953}
\drawhypgeodesic{0.61389}{-0.03953}{0.66158}{0.35281}
\drawhypgeodesic{0.66158}{0.35281}{0.31845}{0.62996}
\drawhypgeodesic{0.31845}{0.62996}{-0.28126}{0.40236}
\drawhypgeodesic{-0.28126}{0.40236}{-0.70212}{0.09548}
\drawhypgeodesic{-0.70212}{0.09548}{-0.64160}{-0.35274}
\drawhypgeodesic{-0.64160}{-0.35274}{-0.25001}{-0.51389}

\filldraw[fill=white] (-0.25001,-0.51389) circle (0.05);      
\end{scope}
\begin{scope}[shift={(9.2,-6.9)}]
\drawhypgeodesic{0.25619}{-0.62827}{0.56840}{-0.42496}
\drawhypgeodesic{0.56840}{-0.42496}{0.68481}{-0.02932}
\drawhypgeodesic{0.68481}{-0.02932}{0.47192}{0.42482}
\drawhypgeodesic{0.47192}{0.42482}{-0.03421}{0.61722}
\drawhypgeodesic{-0.03421}{0.61722}{-0.55475}{0.39569}
\drawhypgeodesic{-0.55475}{0.39569}{-0.73013}{0.03445}
\drawhypgeodesic{-0.73013}{0.03445}{-0.66223}{-0.38965}
\drawhypgeodesic{-0.66223}{-0.38965}{0.25619}{-0.62827}

\filldraw[fill=white] (0.25619,-0.62827) circle (0.05);      
\end{scope}

\begin{scope}[shift={(0,-6.9)}]

\foreach \a in {0,...,7}
\draw
	(45*\a-22.5:0.8) -- (45*\a+22.5:0.8);

\filldraw[fill=white] (22.5:0.8) circle (0.05);

\end{scope}

\end{tikzpicture}
\caption{Geometrical octagons that can tile $3{\bb P}^2$.} 
\label{3P2_n8geom}
\end{figure}

We remark that some settings may produce polygons with self crossing. For example, for tiling 8.6 in Figure \ref{3P2_n8}, if we set 
\begin{align*}
|\bar{0}| &=|\bar{2}|=10a,\;
|\bar{1}|=|\bar{5}|=11a,\;
|\bar{3}|=|\bar{7}|=12a,\;
|\bar{4}|=|\bar{6}|=13a, \\
[0] &=[1]=[2]=[3]=\tfrac{\pi}{2},\;
[4]=x,\;
[5]=y,\;
[6]=\pi-x,\;
[7]=\pi-y,
\end{align*}
then we get the non-simple octagon on the left of Figure \ref{cross}. If we keep the edge length settings and slightly modify the angle settings 
\[
[0]=[2]=\tfrac{2\pi}{3},\;
[1]=[3]=\tfrac{\pi}{3},\;
[4]=x,\;
[5]=y,\;
[6]=\pi-x,\;
[7]=\pi-y,
\]
then we get simple octagon in the middle of Figure \ref{cross}. We may even change the assignment of variable angles, and set 
\[
[0]=[2]=\tfrac{\pi}{2},\;
[3]=[7]=\tfrac{\pi}{3},\;
[1]=x,\;
[2]=y,\;
[5]=\tfrac{4\pi}{3}-x,\;
[6]=\pi-y.
\]
Then we get the simple octagon on the right of Figure \ref{3P2_n8geom}, which is used in Figure \ref{3P2_n8geom}.

\begin{figure}[htp]
\centering
\begin{tikzpicture}[scale=1]

\draw (0,0) circle (1);

\drawhypgeodesic{-0.35463}{-0.55106}{-0.02228}{-0.67182}
\drawhypgeodesic{-0.02228}{-0.67182}{0.34521}{-0.44985}
\drawhypgeodesic{0.34521}{-0.44985}{0.32730}{0.04147}
\drawhypgeodesic{0.32730}{0.04147}{-0.33714}{0.18045}
\drawhypgeodesic{-0.33714}{0.18045}{-0.42830}{0.65978}
\drawhypgeodesic{-0.42830}{0.65978}{-0.43688}{0.82015}
\drawhypgeodesic{-0.43688}{0.82015}{-0.35881}{0.87997}
\drawhypgeodesic{-0.35881}{0.87997}{-0.35463}{-0.55106}

\filldraw[fill=white] (-0.35463,-0.55106) circle (0.05);

\begin{scope}[xshift=3cm]

\draw (0,0) circle (1); 

\drawhypgeodesic{0.03664}{-0.53003}{0.48244}{-0.56215}
\drawhypgeodesic{0.48244}{-0.56215}{0.62115}{-0.12697}
\drawhypgeodesic{0.62115}{-0.12697}{0.65914}{0.29218}
\drawhypgeodesic{0.65914}{0.29218}{0.09292}{0.47172}
\drawhypgeodesic{0.09292}{0.47172}{-0.55381}{0.47731}
\drawhypgeodesic{-0.55381}{0.47731}{-0.78358}{0.21322}
\drawhypgeodesic{-0.78358}{0.21322}{-0.55490}{-0.23528}
\drawhypgeodesic{-0.55490}{-0.23528}{0.03664}{-0.53003}

\filldraw[fill=white] (0.03664,-0.53003) circle (0.05);   
   
\end{scope}

\begin{scope}[xshift=6cm]

\draw (0,0) circle (1); 

\drawhypgeodesic{-0.30537}{-0.57115}{0.19331}{-0.39336}
\drawhypgeodesic{0.19331}{-0.39336}{0.66223}{-0.04977}
\drawhypgeodesic{0.66223}{-0.04977}{0.74470}{0.28824}
\drawhypgeodesic{0.74470}{0.28824}{0.38650}{0.53478}
\drawhypgeodesic{0.38650}{0.53478}{-0.28630}{0.44055}
\drawhypgeodesic{-0.28630}{0.44055}{-0.69267}{0.10843}
\drawhypgeodesic{-0.69267}{0.10843}{-0.70240}{-0.35772}
\drawhypgeodesic{-0.70240}{-0.35772}{-0.30537}{-0.57115}

\filldraw[fill=white] (-0.30537,-0.57115) circle (0.05); 
     
\end{scope}

\end{tikzpicture}
\caption{Non-simple and simple octagons.} 
\label{cross}
\end{figure}

For single tile tilings of orientable surfaces, the combinatorial part of the tiling is uniquely determined by the planar diagram. Therefore choosing different edge lengths for different edge pairs is already enough. Figures \ref{2T2_n8geom} gives the geometric realisations for the octagons and $18$-gons  in Figure \ref{2T2_n8}. The point $\circ$ indicates the vertex of the octagon on the lower left of Figure \ref{3P2_n8geom} (compared with tilings $8.*$ in Figure \ref{2T2_n8}), and the vertex of the 18-gon on the right of Figure \ref{2T2_n8geom} (compared with tilings $18.*$ in Figure \ref{2T2_n16}).

\begin{figure}[htp]
\centering
\begin{tikzpicture}[scale=1]

\foreach \x in {0,...,5}
\foreach \y in {0,1}
\draw (2.3*\x,-2.3*\y) circle (1); 

\foreach \x in {1,...,4}
\node[gray] at (-2.3+2.3*\x, 0) { 8.\x};

\foreach \x in {1,2}
\node[gray] at (2.3*3+2.3*\x, 0) { 18.\x};

\foreach \x in {3,...,8}
\node[gray] at (-2.3*3+2.3*\x, -2.3) { 18.\x};

\foreach \a in {0,...,17}
\draw[shift={(2.3*6,-1.15)}]
	(20*\a:0.9) -- (20*\a+20:0.9);

\filldraw[fill=white] (2.3*6+0.9,-1.15) circle (0.05);
	

\drawhypgeodesic{0.19194}{-0.91487}{0.36630}{-0.89527}
\drawhypgeodesic{0.36630}{-0.89527}{0.48668}{-0.84369}
\drawhypgeodesic{0.48668}{-0.84369}{0.57602}{-0.77776}
\drawhypgeodesic{0.57602}{-0.77776}{0.70339}{-0.59871}
\drawhypgeodesic{0.70339}{-0.59871}{0.65025}{0.15635}
\drawhypgeodesic{0.65025}{0.15635}{-0.47460}{0.69296}
\drawhypgeodesic{-0.47460}{0.69296}{-0.48700}{-0.41518}
\drawhypgeodesic{-0.48700}{-0.41518}{0.19194}{-0.91487}

\filldraw[fill=white] (0.19194,-0.91487) circle (0.05);


\begin{scope}[shift={(2.3,0)}]

\drawhypgeodesic{0.37036}{-0.84775}{0.55976}{-0.78815}
\drawhypgeodesic{0.55976}{-0.78815}{0.68145}{-0.46100}
\drawhypgeodesic{0.68145}{-0.46100}{0.88578}{0.03611}
\drawhypgeodesic{0.88578}{0.03611}{0.80237}{0.48326}
\drawhypgeodesic{0.80237}{0.48326}{0.69317}{0.67232}
\drawhypgeodesic{0.69317}{0.67232}{0.21022}{0.78486}
\drawhypgeodesic{0.21022}{0.78486}{-0.48933}{-0.33305}
\drawhypgeodesic{-0.48933}{-0.33305}{0.37036}{-0.84775}
 
\filldraw[fill=white] (0.37036,-0.84775) circle (0.05);
 
\end{scope}


\begin{scope}[shift={(4.6,0)}]

\drawhypgeodesic{0.38153}{-0.81722}{0.61853}{-0.71283}
\drawhypgeodesic{0.61853}{-0.71283}{0.75478}{-0.60802}
\drawhypgeodesic{0.75478}{-0.60802}{0.81909}{-0.51078}
\drawhypgeodesic{0.81909}{-0.51078}{0.88186}{-0.27306}
\drawhypgeodesic{0.88186}{-0.27306}{0.78556}{0.30507}
\drawhypgeodesic{0.78556}{0.30507}{0.20488}{0.74152}
\drawhypgeodesic{0.20488}{0.74152}{-0.48292}{-0.32733}
\drawhypgeodesic{-0.48292}{-0.32733}{0.38153}{-0.81722}

\filldraw[fill=white] (0.38153,-0.81722) circle (0.05);
  
\end{scope}


\begin{scope}[shift={(6.9,0)}]

\drawhypgeodesic{0.38249}{-0.81548}{0.62155}{-0.70843}
\drawhypgeodesic{0.62155}{-0.70843}{0.77899}{-0.52156}
\drawhypgeodesic{0.77899}{-0.52156}{0.89328}{-0.20034}
\drawhypgeodesic{0.89328}{-0.20034}{0.88217}{0.25078}
\drawhypgeodesic{0.88217}{0.25078}{0.73822}{0.52205}
\drawhypgeodesic{0.73822}{0.52205}{0.20549}{0.73919}
\drawhypgeodesic{0.20549}{0.73919}{-0.48229}{-0.32703}
\drawhypgeodesic{-0.48229}{-0.32703}{0.38249}{-0.81548}
 
\filldraw[fill=white] (0.38249,-0.81548) circle (0.05);
  
\end{scope}


\begin{scope}[shift={(9.2,0)}]

\drawhypgeodesic{0.06818}{-0.76202}{0.25056}{-0.75428}
\drawhypgeodesic{0.25056}{-0.75428}{0.41531}{-0.68405}
\drawhypgeodesic{0.41531}{-0.68405}{0.55822}{-0.55810}
\drawhypgeodesic{0.55822}{-0.55810}{0.64755}{-0.39370}
\drawhypgeodesic{0.64755}{-0.39370}{0.67587}{-0.13885}
\drawhypgeodesic{0.67587}{-0.13885}{0.57942}{0.17443}
\drawhypgeodesic{0.57942}{0.17443}{0.30675}{0.45984}
\drawhypgeodesic{0.30675}{0.45984}{0.09728}{0.73621}
\drawhypgeodesic{0.09728}{0.73621}{-0.08374}{0.83590}
\drawhypgeodesic{-0.08374}{0.83590}{-0.23016}{0.87314}
\drawhypgeodesic{-0.23016}{0.87314}{-0.33015}{0.78248}
\drawhypgeodesic{-0.33015}{0.78248}{-0.45090}{0.60382}
\drawhypgeodesic{-0.45090}{0.60382}{-0.60373}{0.33655}
\drawhypgeodesic{-0.60373}{0.33655}{-0.65510}{0.04510}
\drawhypgeodesic{-0.65510}{0.04510}{-0.58240}{-0.32687}
\drawhypgeodesic{-0.58240}{-0.32687}{-0.43946}{-0.54567}
\drawhypgeodesic{-0.43946}{-0.54567}{-0.22349}{-0.68390}
\drawhypgeodesic{-0.22349}{-0.68390}{0.06818}{-0.76202}

\filldraw[fill=white] (0.06818,-0.76202) circle (0.05);

\end{scope}


\begin{scope}[shift={(2.3*5,0)}]

\drawhypgeodesic{0.08224}{-0.62178}{0.34979}{-0.60542}
\drawhypgeodesic{0.34979}{-0.60542}{0.55942}{-0.48698}
\drawhypgeodesic{0.55942}{-0.48698}{0.70077}{-0.31182}
\drawhypgeodesic{0.70077}{-0.31182}{0.76242}{-0.12571}
\drawhypgeodesic{0.76242}{-0.12571}{0.75838}{0.09685}
\drawhypgeodesic{0.75838}{0.09685}{0.68460}{0.32319}
\drawhypgeodesic{0.68460}{0.32319}{0.55298}{0.49830}
\drawhypgeodesic{0.55298}{0.49830}{0.34592}{0.62696}
\drawhypgeodesic{0.34592}{0.62696}{0.06304}{0.68080}
\drawhypgeodesic{0.06304}{0.68080}{-0.24884}{0.62482}
\drawhypgeodesic{-0.24884}{0.62482}{-0.52598}{0.46165}
\drawhypgeodesic{-0.52598}{0.46165}{-0.71082}{0.22546}
\drawhypgeodesic{-0.71082}{0.22546}{-0.79896}{0.00754}
\drawhypgeodesic{-0.79896}{0.00754}{-0.81908}{-0.17053}
\drawhypgeodesic{-0.81908}{-0.17053}{-0.80620}{-0.34243}
\drawhypgeodesic{-0.80620}{-0.34243}{-0.64209}{-0.43042}
\drawhypgeodesic{-0.64209}{-0.43042}{-0.30759}{-0.45048}
\drawhypgeodesic{-0.30759}{-0.45048}{0.08224}{-0.62178}

\filldraw[fill=white] (0.08224,-0.62178) circle (0.05);
 
\end{scope}


\begin{scope}[shift={(0,-2.3)}]

\drawhypgeodesic{0.11211}{-0.64382}{0.36124}{-0.62504}
\drawhypgeodesic{0.36124}{-0.62504}{0.55683}{-0.51532}
\drawhypgeodesic{0.55683}{-0.51532}{0.69294}{-0.35234}
\drawhypgeodesic{0.69294}{-0.35234}{0.75655}{-0.17424}
\drawhypgeodesic{0.75655}{-0.17424}{0.75689}{0.04721}
\drawhypgeodesic{0.75689}{0.04721}{0.68143}{0.28193}
\drawhypgeodesic{0.68143}{0.28193}{0.52406}{0.50158}
\drawhypgeodesic{0.52406}{0.50158}{0.32508}{0.63546}
\drawhypgeodesic{0.32508}{0.63546}{0.06615}{0.68438}
\drawhypgeodesic{0.06615}{0.68438}{-0.24206}{0.62075}
\drawhypgeodesic{-0.24206}{0.62075}{-0.51858}{0.44331}
\drawhypgeodesic{-0.51858}{0.44331}{-0.72987}{0.24210}
\drawhypgeodesic{-0.72987}{0.24210}{-0.82263}{0.07130}
\drawhypgeodesic{-0.82263}{0.07130}{-0.87018}{-0.09813}
\drawhypgeodesic{-0.87018}{-0.09813}{-0.77846}{-0.23983}
\drawhypgeodesic{-0.77846}{-0.23983}{-0.60733}{-0.38661}
\drawhypgeodesic{-0.60733}{-0.38661}{-0.26415}{-0.49272}
\drawhypgeodesic{-0.26415}{-0.49272}{0.11211}{-0.64382}
  
\filldraw[fill=white] (0.11211,-0.64382) circle (0.05);
 
\end{scope}


\begin{scope}[shift={(2.3,-2.3)}]

\drawhypgeodesic{-0.80376}{-0.49875}{-0.71891}{-0.53936}
\drawhypgeodesic{-0.71891}{-0.53936}{-0.54244}{-0.58944}
\drawhypgeodesic{-0.54244}{-0.58944}{-0.19845}{-0.59395}
\drawhypgeodesic{-0.19845}{-0.59395}{0.29062}{-0.40251}
\drawhypgeodesic{0.29062}{-0.40251}{0.59430}{-0.05619}
\drawhypgeodesic{0.59430}{-0.05619}{0.72031}{0.32466}
\drawhypgeodesic{0.72031}{0.32466}{0.74306}{0.53466}
\drawhypgeodesic{0.74306}{0.53466}{0.74166}{0.61447}
\drawhypgeodesic{0.74166}{0.61447}{0.73941}{0.65349}
\drawhypgeodesic{0.73941}{0.65349}{0.69321}{0.65362}
\drawhypgeodesic{0.69321}{0.65362}{0.56187}{0.60954}
\drawhypgeodesic{0.56187}{0.60954}{0.51740}{0.14855}
\drawhypgeodesic{0.51740}{0.14855}{-0.18479}{0.27531}
\drawhypgeodesic{-0.18479}{0.27531}{-0.62961}{0.02162}
\drawhypgeodesic{-0.62961}{0.02162}{-0.79822}{-0.25934}
\drawhypgeodesic{-0.79822}{-0.25934}{-0.84522}{-0.41115}
\drawhypgeodesic{-0.84522}{-0.41115}{-0.85981}{-0.47195}
\drawhypgeodesic{-0.85981}{-0.47195}{-0.80376}{-0.49875}

\filldraw[fill=white] (-0.80376,-0.49875) circle (0.05);
 
\end{scope}


\begin{scope}[shift={(2.3*2,-2.3)}]

\drawhypgeodesic{0.03975}{-0.76776}{0.23027}{-0.75431}
\drawhypgeodesic{0.23027}{-0.75431}{0.39138}{-0.67861}
\drawhypgeodesic{0.39138}{-0.67861}{0.54906}{-0.52312}
\drawhypgeodesic{0.54906}{-0.52312}{0.64214}{-0.32375}
\drawhypgeodesic{0.64214}{-0.32375}{0.66521}{0.00242}
\drawhypgeodesic{0.66521}{0.00242}{0.57251}{0.33075}
\drawhypgeodesic{0.57251}{0.33075}{0.38587}{0.57467}
\drawhypgeodesic{0.38587}{0.57467}{0.18893}{0.72865}
\drawhypgeodesic{0.18893}{0.72865}{0.00360}{0.77936}
\drawhypgeodesic{0.00360}{0.77936}{-0.20543}{0.76369}
\drawhypgeodesic{-0.20543}{0.76369}{-0.38790}{0.69395}
\drawhypgeodesic{-0.38790}{0.69395}{-0.52299}{0.54494}
\drawhypgeodesic{-0.52299}{0.54494}{-0.64294}{0.27517}
\drawhypgeodesic{-0.64294}{0.27517}{-0.67333}{-0.03378}
\drawhypgeodesic{-0.67333}{-0.03378}{-0.59128}{-0.35012}
\drawhypgeodesic{-0.59128}{-0.35012}{-0.43724}{-0.55845}
\drawhypgeodesic{-0.43724}{-0.55845}{-0.20762}{-0.70371}
\drawhypgeodesic{-0.20762}{-0.70371}{0.03975}{-0.76776}
  
\filldraw[fill=white] (0.03975,-0.76776) circle (0.05);
 
\end{scope}


\begin{scope}[shift={(2.3*3,-2.3)}]

\drawhypgeodesic{-0.29997}{-0.69013}{-0.03488}{-0.71046}
\drawhypgeodesic{-0.03488}{-0.71046}{0.22678}{-0.63177}
\drawhypgeodesic{0.22678}{-0.63177}{0.51339}{-0.40565}
\drawhypgeodesic{0.51339}{-0.40565}{0.67516}{-0.13542}
\drawhypgeodesic{0.67516}{-0.13542}{0.73326}{0.15871}
\drawhypgeodesic{0.73326}{0.15871}{0.70924}{0.37432}
\drawhypgeodesic{0.70924}{0.37432}{0.64169}{0.51256}
\drawhypgeodesic{0.64169}{0.51256}{0.50697}{0.61832}
\drawhypgeodesic{0.50697}{0.61832}{0.30523}{0.69346}
\drawhypgeodesic{0.30523}{0.69346}{0.02980}{0.70876}
\drawhypgeodesic{0.02980}{0.70876}{-0.23258}{0.62673}
\drawhypgeodesic{-0.23258}{0.62673}{-0.51787}{0.39500}
\drawhypgeodesic{-0.51787}{0.39500}{-0.68449}{0.11275}
\drawhypgeodesic{-0.68449}{0.11275}{-0.74281}{-0.17269}
\drawhypgeodesic{-0.74281}{-0.17269}{-0.71763}{-0.35522}
\drawhypgeodesic{-0.71763}{-0.35522}{-0.63126}{-0.49233}
\drawhypgeodesic{-0.63126}{-0.49233}{-0.48005}{-0.60694}
\drawhypgeodesic{-0.48005}{-0.60694}{-0.29997}{-0.69013}
  
\filldraw[fill=white] (-0.29997,-0.69013) circle (0.05);
 
\end{scope}


\begin{scope}[shift={(2.3*4,-2.3)}]

\drawhypgeodesic{0.12193}{-0.65591}{0.37219}{-0.62122}
\drawhypgeodesic{0.37219}{-0.62122}{0.58368}{-0.49783}
\drawhypgeodesic{0.58368}{-0.49783}{0.70226}{-0.35327}
\drawhypgeodesic{0.70226}{-0.35327}{0.77151}{-0.16297}
\drawhypgeodesic{0.77151}{-0.16297}{0.77666}{0.02730}
\drawhypgeodesic{0.77666}{0.02730}{0.70813}{0.23114}
\drawhypgeodesic{0.70813}{0.23114}{0.53524}{0.45493}
\drawhypgeodesic{0.53524}{0.45493}{0.28970}{0.61016}
\drawhypgeodesic{0.28970}{0.61016}{-0.00303}{0.65645}
\drawhypgeodesic{-0.00303}{0.65645}{-0.34970}{0.57557}
\drawhypgeodesic{-0.34970}{0.57557}{-0.61611}{0.44777}
\drawhypgeodesic{-0.61611}{0.44777}{-0.75238}{0.31163}
\drawhypgeodesic{-0.75238}{0.31163}{-0.83672}{0.15860}
\drawhypgeodesic{-0.83672}{0.15860}{-0.81032}{-0.03534}
\drawhypgeodesic{-0.81032}{-0.03534}{-0.74087}{-0.19816}
\drawhypgeodesic{-0.74087}{-0.19816}{-0.54458}{-0.40021}
\drawhypgeodesic{-0.54458}{-0.40021}{-0.20759}{-0.54863}
\drawhypgeodesic{-0.20759}{-0.54863}{0.12193}{-0.65591}
  
\filldraw[fill=white] (0.12193,-0.65591) circle (0.05);
 
\end{scope}


\begin{scope}[shift={(2.3*5,-2.3)}]

\drawhypgeodesic{0.09048}{-0.74581}{0.28288}{-0.72933}
\drawhypgeodesic{0.28288}{-0.72933}{0.45129}{-0.64634}
\drawhypgeodesic{0.45129}{-0.64634}{0.58975}{-0.50509}
\drawhypgeodesic{0.58975}{-0.50509}{0.68532}{-0.30585}
\drawhypgeodesic{0.68532}{-0.30585}{0.71374}{-0.05300}
\drawhypgeodesic{0.71374}{-0.05300}{0.64272}{0.18067}
\drawhypgeodesic{0.64272}{0.18067}{0.40722}{0.43023}
\drawhypgeodesic{0.40722}{0.43023}{0.16222}{0.67292}
\drawhypgeodesic{0.16222}{0.67292}{-0.05604}{0.75205}
\drawhypgeodesic{-0.05604}{0.75205}{-0.27484}{0.74769}
\drawhypgeodesic{-0.27484}{0.74769}{-0.47003}{0.69438}
\drawhypgeodesic{-0.47003}{0.69438}{-0.58213}{0.54246}
\drawhypgeodesic{-0.58213}{0.54246}{-0.68273}{0.30800}
\drawhypgeodesic{-0.68273}{0.30800}{-0.70652}{0.07120}
\drawhypgeodesic{-0.70652}{0.07120}{-0.62168}{-0.22978}
\drawhypgeodesic{-0.62168}{-0.22978}{-0.43205}{-0.51331}
\drawhypgeodesic{-0.43205}{-0.51331}{-0.19962}{-0.67109}
\drawhypgeodesic{-0.19962}{-0.67109}{0.09048}{-0.74581}
   
\filldraw[fill=white] (0.09048,-0.74581) circle (0.05);
 
\end{scope}

\end{tikzpicture}
\caption{Geometrical octagons and $18$-gons that can tile $2{\bb T}^2$.} 
\label{2T2_n8geom}
\end{figure}

We verified the geometrical realisation for all single tile tilings of $2{\bb T}^2$ and $3{\bb P}^2$, in Figures \ref{2T2_n8}, \ref{2T2_n16}, and \ref{3P2_n8}.

\end{document}